\documentclass[A4,11pt]{article}
\usepackage{hyperref}
\usepackage{amssymb,amsmath}
\usepackage{color}
\usepackage{lmodern}
\usepackage{bm}
\usepackage{amssymb,color}
\usepackage{titlesec}
\titleformat{\paragraph}[block]{\normalfont\bfseries}{\theparagraph}{1em}{}
\usepackage{mathptmx}  
\usepackage{amsthm}

\def\Ph{\phantom{\int}\!\!\!\!\!}

\def\Rd{\mathbb{R}^{d}}

\def\black#1{{#1}}

\def\wphi{\widetilde{\Phi}}

\def\da{\frac{d}{\alpha}}
\def\C{{\cal C}}

\def\E{{\cal S}(\Rd)}
\def\ese{{\cal S}(\mathbb{R}^{d+1})}

\def\J{\textbf{\textit{X}}}
\def\X{\textbf{\textit{X}}}

\def\text#1{{\rm #1}}
\def\EE{\mathbb{E}}
\def\L{\left}
\def\R{\right}
\def\b{3-\da}

\newtheorem{theorem}{Theorem}[section]  
\newtheorem{corollary}[theorem]{Corollary}

\newtheorem{remark}[theorem]{Remark}

\begin{document}

	\title{Occupation-Time Fluctuations of  an Age-Dependent Branching System Driven by Poisson Immigration\vspace{-1ex}}
\author{Ekaterina T. Kolkovska\thanks{Centro de Investigación en Matemáticas, Guanajuato, Mexico.} \and José Alfredo López-Mimbela\thanks{Centro de Investigación en Matemáticas, Guanajuato, Mexico.} \and Maroussia Slavtchova-Bojkova\thanks{Faculty of Mathematics and Informatics, Sofia University “St. Kliment Ohridski”, and the Institute of Mathematics and Informatics, Bulgarian Academy of Sciences, Sofia, Bulgaria.}}
\date{}
	\maketitle

\begin{abstract}
We study the occupation-time fluctuations of a critical age-dependent branching particle system with immigration in $\mathbb{R}^d$. Immigrants arrive according to a homogeneous Poisson random measure in space and time. Each particle moves independently according to a symmetric $\alpha$-stable process and, at the end of its lifetime, either dies or splits into two offspring with equal probability. The lifetime distribution is allowed to have either finite mean or a heavy tail of index $\gamma\in(0,1]$.
We investigate the asymptotic behavior of the centered occupation-time process under a suitable space-time scaling. Assuming
$
\alpha<d<(1+\gamma)\alpha,
$
we prove that the rescaled occupation-time fluctuations weakly converge as processes with values in the space of tempered distributions to a centered Gaussian process with an explicitly identified covariance structure. The normalization and the covariance depend on both the stability index $\alpha$ and the tail exponent $\gamma$.
The limiting process is self-similar, possesses long-range dependence, and is neither Markovian nor a semimartingale. In contrast with the corresponding age-dependent branching system without immigration, the contribution of the initial population vanishes in the limit, and the asymptotic fluctuations are entirely determined by the immigration mechanism. When $\gamma=1$, our results recover the covariance structure previously obtained for branching systems with immigration and finite-mean lifetimes.
The proofs rely on the space-time random field approach, Fourier analytic techniques, and asymptotic properties of renewal functions associated with the lifetime distribution.
\end{abstract}

\section{Introduction}

In this paper we investigate the occupation-time fluctuations of a critical age-dependent branching particle system with immigration. Branching particle systems constitute a fundamental class of stochastic spatial models describing the evolution of spatially distributed populations subject to reproduction, migration, and random environmental effects. Their large-scale asymptotic behavior has been extensively studied over the last decades and has revealed deep connections between branching mechanisms, particle motion, random fields, and self-similar Gaussian processes.

A central problem in this theory is to understand how the interplay between branching and spatial dispersion determines the fluctuations of occupation times. In many classical models, suitably rescaled occupation-time processes converge to Gaussian limits whose temporal structure reflects both the genealogical structure of the population and the recurrence properties of the underlying motion. Depending on the dimension and the branching mechanism, the limiting processes may exhibit long-range dependence, self-similarity, and other non-Markovian features.

The situation becomes substantially richer when one allows particles to have general lifetime distributions. In age-dependent branching systems, the lifetime law introduces an additional source of memory that may significantly alter the asymptotic behavior of the population. In particular, heavy-tailed lifetimes generate long temporal correlations and lead to fluctuation limits that differ markedly from those arising in Markov branching systems.

Occupation-time fluctuations of branching systems with exponentially distributed lifetimes have been extensively studied by Cox and Griffeath \cite{CG}, Méléard and Roelly \cite{SS}, and subsequently by Bojdecki, Gorostiza, and Talarczyk \cite{BGT,BGT1,BGT2}, with particular emphasis on the rescaled occupation-time process. These works established a general framework showing that a broad class of Gaussian processes emerges as scaling limits, while also highlighting the crucial role of the spatial dimension and the underlying particle motion. 
The model with particle immigration and exponentially distributed lifetimes was studied in \cite{LER}. The case of critical binary age-dependent branching particle systems without immigration and with lifetimes having a power-law tail of order $\gamma^{-1}$ (with $\gamma \le 1$) was investigated in \cite{ALEA}. We further note that immigration mechanisms in branching and measure-valued processes have been systematically studied in the work of Li; see, for instance, \cite{Li1992,LiBook}.

Under the assumption that the lifetime distribution satisfies
\[
1-F(t)\sim \frac{1}{\Gamma(1-\gamma)}\, t^{-\gamma},
\qquad t\to\infty,
\]
for some $0<\gamma\le 1$, where $\Gamma$ denotes the Gamma function, functional limit theorems for the occupation-time fluctuation process were established in \cite{ALEA}. In particular, it was shown that the temporal structure of the limiting process is described by weighted sub-fractional Brownian motion. These results demonstrate that heavy-tailed lifetimes substantially alter both the scaling exponents and the critical dimensions of the model.

The present work may be viewed as a natural continuation of the line of research initiated in \cite{ALEA}; see also \cite{LER}. Our main objective is to understand how the introduction of immigration influences the fluctuation behavior in the presence of heavy-tailed lifetimes. More specifically, we investigate whether the long-range temporal dependence induced by heavy-tailed lifetimes persists under immigration, and how this additional mechanism impacts the scaling exponents and the structure of the limiting process.

Before introducing the spatial model considered in this paper, we briefly recall some relevant developments on branching processes with immigration in the non-spatial setting. Early contributions include the work of Sevastyanov \cite{Sevastianov}, who introduced branching models with immigration driven by Poisson random measures. Many foundational results in this area were later surveyed by Vatutin and Zubkov \cite{7.,8.}. Subsequently, models with time-dependent immigration rates were investigated by Mitov, Yanev, Hyrien, and Slavtchova-Bojkova in \cite{1.,2.,3.}. More recent contributions by Barczy et al. \cite{9.,10.}, González {\it et al.} \cite{11.}, and Li {\it et al.} \cite{12.} have further broadened the theoretical framework. See also Rahimov’s recent survey \cite{17.} for a comprehensive overview.

We consider a critical age-dependent branching particle system evolving in \(\mathbb{R}^d\), in which immigrants arrive according to a homogeneous Poisson random measure on space-time. Once born, particles move independently according to a symmetric \(\alpha\)-stable Lévy process and, at the end of their lifetimes, either die or split into two offspring with equal probability. The lifetime distribution is allowed to have either finite mean or a heavy tail with index \(\gamma\in(0,1]\).

The introduction of immigration fundamentally changes the asymptotic behavior of the system. While age-dependent branching systems without immigration retain a strong memory of the initial population, the present work shows that, after suitable normalization, the contribution of the initial population becomes asymptotically negligible. The limiting fluctuations are generated entirely by the immigration mechanism. This phenomenon is not directly apparent from the particle dynamics and only becomes visible at the level of the fluctuation limits.

Our main result is a functional central limit theorem for the occupation-time fluctuation process in the intermediate-dimensional regime
\[
\alpha<d<(1+\gamma)\alpha .
\]
After an appropriate normalization, the fluctuation process converges weakly in
\(
C\big([0,\Upsilon],\mathcal S'(\mathbb R^d)\big)
\)
to a centered Gaussian process with an explicitly identified covariance structure. The covariance depends on both the stability index \(\alpha\) and the lifetime parameter \(\gamma\), and reduces, when \(\gamma=1\), to the covariance previously obtained for branching systems with immigration and exponentially-distributed lifetimes \cite{LER}.

The limiting Gaussian process exhibits several remarkable properties. It is self-similar, possesses long-range dependence, and is neither Markovian nor a semimartingale. These features indicate that the memory introduced by the lifetime distribution survives in the scaling limit despite the continuous inflow of immigrants. At the same time, however, the immigration mechanism completely dominates the contribution of the initial population, producing fluctuation behavior that differs substantially from that of the corresponding system without immigration.

Our approach is based on the space-time random field method introduced by Bojdecki, Gorostiza, and Ramaswamy \cite{BGR}. This framework reduces the problem to the analysis of Laplace functionals in spaces of tempered distributions. The proof combines Fourier analytic techniques with asymptotic properties of renewal functions associated with the lifetime distribution and delicate estimates involving stable semigroups.

The paper is organized as follows. Section~\ref{Background} introduces the model, the space-time random field framework and collects the moment identities needed in the sequel. Section~\ref{section.results} contains the statement of the main results. The convergence of the covariance functionals and the proof of the functional limit theorem are established in Section~\ref{CovaConver}. Finally, Section~\ref{ProofProperties} is devoted to the study of the structural properties of the limiting Gaussian process, including self-similarity, Hölder regularity, long-range dependence, and non-Markovianity.

\section{Background}\label{Background}

Throughout this section and henceforth, $\mathcal{S}(\mathbb{R}^d)$ denotes the Schwartz space of smooth, rapidly decreasing functions on $\mathbb{R}^d$, and $\mathcal{S}'(\mathbb{R}^d)$ its strong dual.

\subsection{Description of the Branching System}

We assume that immigrant particles enter the population according to a
Poisson random measure on $\mathbb{R}^d\times\mathbb{R}_+$ with intensity
measure
$
\lambda\times\lambda_{[0,\infty)}.
$
Here $\lambda$ and $\lambda_{[0,\infty)}$ denote the Lebesgue measures on
$\mathbb{R}^d$ and $\mathbb{R}_+$, respectively.
Each particle evolves independently according to the following rules:
\begin{enumerate}
\item The offspring number $\xi$ has probability generating function
\[
h(s):=\mathbb{E}s^\xi=\frac12+\frac{s^2}{2},
\qquad |s|\le1,
\]
and offspring particles are born at the death location of the parent.

\item The lifetime distribution is characterized by its renewal function
${\cal U}$, which is assumed to be absolutely continuous.

\item Particle motion is governed by a symmetric $\alpha$-stable process
$\{W_t^\alpha,\ t\ge0\}$ with transition densities
$\{p_t^\alpha,\ t>0\}$ and semigroup $\{{\cal T}_t,\ t\ge0\}$ given by
\[
{\cal T}_t\varphi(x)
=
\int_{\mathbb{R}^d}p_t^\alpha(x-y)\varphi(y)\,dy,
\qquad x\in\mathbb{R}^d,\quad \varphi\in\E.
\]
\end{enumerate}

We denote by $N_t$ the random counting measure on $\mathbb{R}^d$
describing the particle configuration at time $t\ge0$.

\subsection{Moment Identities}

The characteristic functional of the process
$N=\{N_t,\ t\ge0\}$ admits the following representation: for
$0\le t_1\le\cdots\le t_m$ and $\varphi_j\in\E$,
$j=1,\ldots,m$,
\begin{equation}\label{ChF}
\mathbb{E}
\exp\left\{
i\sum_{j=1}^m \langle N_{t_j},\varphi_j\rangle
\right\}
=
\exp\left\{
\int_{\mathbb{R}^d\times[0,\infty)}
\left(
\mathbb{E}
\exp\left\{
i\sum_{j=1}^m
\langle N_{t_j-s}(x),\varphi_j\rangle
\right\}
-1
\right)
\,dx\,ds
\right\}.
\end{equation}
Here $N_r(x)$ denotes the random measure generated at time $r\ge0$
by a branching system starting from a single particle located at
$x\in\mathbb{R}^d$, and we set $N_r(x)=0$ for $r<0$.
Differentiating \eqref{ChF} yields
\begin{align*}
\mathbb{E}\langle N_t,\varphi\rangle
&=
\int_{\mathbb{R}^d}\int_0^t
\mathbb{E}\langle N_{t-r}(x),\varphi\rangle\,dr\,dx,
\\
\mathrm{Cov}\Big(
\langle N_s,\varphi\rangle,
\langle N_t,\psi\rangle
\Big)
&=
\int_{\mathbb{R}^d}\int_0^s
\mathbb{E}\Big[
\langle N_{s-r}(x),\varphi\rangle
\langle N_{t-r}(x),\psi\rangle
\Big]\,dr\,dx,
\end{align*}
for $t\ge s\ge0$ and $\varphi,\psi\in\E$.
By Lemma~3.2 in \cite{ESAIM},
\begin{align*}
\mathbb{E}\langle N_t(x),\varphi\rangle
&=
{\cal T}_t\varphi(x),
\\[4pt]
\mathbb{E}\Big[
\langle N_s(x),\varphi\rangle
\langle N_t(x),\psi\rangle
\Big]
&=
{\cal T}_s\big(\varphi\,{\cal T}_{t-s}\psi\big)(x)
+
\int_0^s
{\cal T}_r\Big[
({\cal T}_{s-r}\varphi)
({\cal T}_{t-r}\psi)
\Big](x)\,d{\cal U}(r).
\end{align*}

Using the invariance of $\lambda$ under $\{{\cal T}_t,\ t\ge0\}$, we obtain
\begin{align}
\nonumber
\mathbb{E}\langle N_t,\varphi\rangle
&=
t\int_{\mathbb{R}^d}\varphi(x)\,dx,
\\
\label{CCOO}
\mathrm{Cov}\Big(
\langle N_s,\varphi\rangle,
\langle N_t,\psi\rangle
\Big)
&=
s\int_{\mathbb{R}^d}
\varphi(x)\,{\cal T}_{t-s}\psi(x)\,dx
\\
\nonumber
&\quad+
\int_0^s
\int_0^{s-u}
\int_{\mathbb{R}^d}
({\cal T}_{s-u-r}\varphi)(x)
({\cal T}_{t-u-r}\psi)(x)
\,dx\,d{\cal U}(r)\,du.
\end{align}

We use the Fourier transform convention
\[
\widehat{\varphi}(z)
=
\int_{\mathbb{R}^d}e^{iz\cdot x}\varphi(x)\,dx,
\qquad z\in\mathbb{R}^d.
\]
Then
\[
\widehat{{\cal T}_r\varphi}(z)
=
e^{-r|z|^\alpha}\widehat{\varphi}(z),
\qquad z\in\mathbb{R}^d,\quad r\ge0.
\]
Thus, by the Parseval--Plancherel identity, formula \eqref{CCOO} becomes
\begin{equation}
\label{Plancherel}
\mathrm{Cov}\Big(
\langle N_s,\varphi\rangle,
\langle N_t,\psi\rangle
\Big)
=
\frac{1}{(2\pi)^d}
\int_{\mathbb{R}^d}
\widehat{\varphi}(z)\overline{\widehat{\psi}(z)}
\Bigg(
s e^{-(t-s)|z|^\alpha}
+
\int_0^s
\int_0^{s-u}
e^{-(t+s-2u-2r)|z|^\alpha}
\,d{\cal U}(r)\,du
\Bigg)\,dz.
\end{equation}

\subsection{Occupation-Time Fluctuations and the Space-Time Random Field Approach}

For each $T>0$, define the \emph{occupation-time fluctuation process}
$\X_T=\{\X_T(t),\ t\ge0\}$ by
\begin{equation}\label{DefX}
\langle \X_T(t),\varphi\rangle
=
\frac{1}{F_T}
\left(
\int_0^{Tt}\langle N_s,\varphi\rangle\,ds
-
\mathbb{E}\int_0^{Tt}\langle N_s,\varphi\rangle\,ds
\right),
\qquad \varphi\in\E,
\end{equation}
where $F_T$ is a normalizing factor chosen so that $\X_T$
converges weakly to a non-trivial limit as $T\to\infty$.
Observe that, because of the centering, $\X_T$ is no longer
measure-valued. Instead, it is viewed as a generalized process taking
values in $\mathcal{S}'(\mathbb{R}^d)$.
We now recall a useful criterion for proving weak convergence of generalized
processes. Fix $\Upsilon>0$ and consider the space
$
C([0,\Upsilon],\mathcal{S}'(\mathbb{R}^d)),
$
endowed with the topology of uniform convergence. For each
$
x\in C([0,\Upsilon],\mathcal{S}'(\mathbb{R}^d)),
$
define $\widetilde{x}\in\mathcal{S}'(\mathbb{R}^{d+1})$ by
\begin{equation}\label{DefDelFuncional}
\langle \widetilde{x},\Phi\rangle
=
\int_0^\Upsilon
\langle x(t),\Phi(\cdot,t)\rangle\,dt,
\qquad
\Phi\in\mathcal{S}(\mathbb{R}^{d+1}).
\end{equation}
This mapping is continuous. Therefore, every process $X$ with sample paths in
$
C([0,\Upsilon],\mathcal{S}'(\mathbb{R}^d))
$
induces an $\mathcal{S}'(\mathbb{R}^{d+1})$-valued random variable
$\widetilde X$, called its \emph{space-time random field}. Moreover, the law
of $X$ is uniquely determined by the law of $\widetilde X$; see \cite{BGR}.
Weak convergence is denoted by $\Rightarrow$. The following criterion is a
particular case of \cite[Theorem 4.3]{BGR}.

\begin{theorem}\label{RFA}
Let $\{X_n,\ n\ge1\}$ be a sequence of processes with paths in
$
C([0,\Upsilon],\mathcal{S}'(\mathbb{R}^d)).
$
Assume that:
\begin{enumerate}
\item $\{X_n,\ n\ge1\}$ is tight;
\item $\{\widetilde X_n,\ n\ge1\}$ converges weakly in
$\mathcal{S}'(\mathbb{R}^{d+1})$.
\end{enumerate}
Then there exists an $\mathcal{S}'(\mathbb{R}^d)$-valued process $X$ such that
$
X_n\Rightarrow X.
$
Moreover, the law of $X$ is uniquely determined by the limit law of
$\widetilde X_n$.
\end{theorem}

\section{Main results}\label{section.results}
Recall that, for each $T>0$, the occupation-time fluctuation process
$\X_T=\{\X_T(t),\,t\ge0\}$ is defined by
\eqref{DefX}.
\begin{theorem}[Convergence in path space]\label{Covariance convergence}
Let $F$ be an absolutely continuous lifetime distribution such that
\[
F(0)=0,\quad F(x)<1 \ \textit{\ for \ all } \ x\ge 0,\quad \textit{and}\quad
1-F(t)\sim\frac{1}{t^{\gamma}\Gamma(1-\gamma)}\quad\mbox{with}\quad\gamma<1.
\]
Let $0<\alpha<d<(1+\gamma)\alpha$, and $F_T^2 = T^{3+\gamma-\frac{d}{\alpha}}.$
Then $\X_T$ converges weakly to $\X$ in $C([0,\Upsilon],{\cal S}'(\mathbb{R}^d))$ as $T\to\infty$ for any $\Upsilon>0$, where $\X:=\{\X(t),\ t\ge 0\}$ is a centered Gaussian process with covariance functional given by
$$
\mathrm{Cov}\big(\langle \X(s),\varphi\rangle,\langle \X(t),\psi\rangle\big)
=
\frac{\gamma\Gamma(1/\alpha)\langle\lambda,\varphi\rangle\langle\lambda,\psi\rangle}
{2(\frac{d}{\alpha}-1)(2-\frac{d}{\alpha})(3-\frac{d}{\alpha})
\pi\alpha\Gamma(1+\gamma)} C(s,t),
$$
where
\begin{align}
\nonumber
C(s,t)
&=
\int_0^s w^{\gamma-1}\Bigg[
3(s-w)^{\b} + (t-w)^{\b}
+\frac{\b}{\gamma}w\Big(
(t-w)^{2-\frac{d}{\alpha}} - (s-w)^{2-\frac{d}{\alpha}}
\Big)
-(t+s-2w)^{\b}
\Bigg]\,dw,
\\ \label{F2}
&\qquad 0\le s\le t,\quad \varphi,\psi\in\E.
\end{align}
\end{theorem}
\begin{corollary}
Assume that $\gamma=1$ and let
\(
h:=3-d/\alpha.
\)
Then, for $0\le s\le t$,
$$
\mathrm{Cov}\big(\langle \X(s),\varphi\rangle,\langle \X(t),\psi\rangle\big)
=
\frac{\Gamma(1/\alpha)\langle\lambda,\varphi\rangle\langle\lambda,\psi\rangle}
{2(\frac{d}{\alpha}-1)(h-1)h
\pi\alpha}
\L[
\frac{
2s^{h+1}
+
2t^{h+1}
+
\left(h-\frac12\right)(t-s)^{h+1}
-
\frac12(t+s)^{h+1}
}{h+1}
-
t(t-s)^h\R].
$$
\end{corollary}

\begin{remark}{\em 
We emphasize that the branching system considered in this paper starts from an empty initial population. This choice is motivated by the following observation.

Suppose that, in addition to the immigration mechanism, the system is endowed with an independent initial Poisson population on $\Rd$ with homogeneous intensity. The corresponding model without immigration was studied in \cite{ALEA}, where it was shown that, in the regime
\(
\alpha\gamma < d < \alpha(1+\gamma)
\) with $d\neq\alpha$,
the occupation-time fluctuation process converges under the normalization
\(
H_T^2=T^{2+\gamma-d/\alpha}.
\)
By contrast, in the present model with immigration, the correct normalization is
\(
F_T^2=T^{3+\gamma-d/\alpha}.
\)
Since
\[
\frac{H_T^2}{F_T^2}=T^{-1}\longrightarrow 0,
\qquad T\to\infty,
\]
the contribution of the initial population is negligible under the scaling considered here. Consequently, the limiting occupation-time fluctuations are entirely determined by the immigration mechanism.
}
\end{remark}
\begin{theorem}[Properties of the limit process]\label{Properties}$\,$  
  \begin{enumerate}
    \item The limit process $\X$ is self-similar with self-similarity index
$\frac{3+\gamma-d/\alpha}{2}$, that is, for every $a>0$,
\[
\{\X(at),\ t\ge0\}
\stackrel{d}{=}
\left\{a^{\frac{3+\gamma-d/\alpha}{2}}\X(t),\ t\ge0\right\}.
\]

\item For every $0\le s<t$, we have $C(s,t)>0$. In particular, the covariance functional of $\X$ is strictly positive.

\item For every $0\le s\le t$,
\begin{equation}\label{nosemimtgla}
\frac{1}{\gamma}\left(4-2^{3-\frac{d}{\alpha}}\right) s^{\gamma}(t-s)^{3-\frac{d}{\alpha}}
\le
\EE\big(\X(t)-\X(s)\big)^2
\le
\frac{2}{\gamma} (t-s)^{3-\frac{d}{\alpha}}.
\end{equation}
In particular, $\X$ is not a semimartingale.

\item $\X$ has Hölder-continuous sample paths of any order less than
$\frac{3-d/\alpha}{2}$; that is, for every
\(
\varepsilon\in\left(0,\frac{3-d/\alpha}{2}\right)
\)
and every $L>0$, there exists a random variable $K_{\varepsilon,L}\ge0$ such that, almost surely,
\[
|\X(t)-\X(s)|
\le
K_{\varepsilon,L}
|t-s|^{(3-d/\alpha)/2-\varepsilon},
\qquad s,t\in[0,L].
\]

\item $\X$ is not Markovian.

\item $\X$ exhibits long-range dependence  in the sense that, for
$0\le r<v\le s<t$,
\[
\lim_{T\to\infty}
T^{d/\alpha}
\EE\big[(\X(t+T)-\X(s+T))(\X(v)-\X(r))\big]
=
\frac{2-d/\alpha}{\gamma(1+\gamma)}
(t-s)\big(v^{1+\gamma}-r^{1+\gamma}\big).
\]
  \end{enumerate}
\end{theorem}

\section{Proof of Theorem \ref{Covariance convergence}}\label{CovaConver}
	
	We start by showing convergence of the covariance functionals which leads to the expression \eqref{F2}.
		From the definition of the process $\X_T$,
		\begin{align*}
	{\rm Cov}\L(\L\langle\X_T(s),\varphi\R\rangle,\L\langle\X_T(t),\psi\R\rangle\R)		& = \frac{T^2}{F^2_T}\int_0^sdu\int_0^t\,dv\:{\rm Cov}\L(\langle N_{Tu},\varphi\rangle,\langle N_{Tv},\psi\rangle\R)               \\
			& = \frac{T^2}{F^2_T}\L(\int_0^sdu\int_s^tdv+\int_0^sdu\int_u^sdv+\int_0^sdv\int_v^sdu\R)                                     {\rm Cov}\L(\langle N_{Tu},\varphi\rangle,\langle N_{Tv},\psi\rangle\R)                                                         \\
			& =  \frac{T^2}{F^2_T}
			\L(\int_0^sdu\int_s^tdv+2\int_0^sdu\int_u^sdv\R)
			{\rm Cov}\L(\langle N_{Tu},\varphi\rangle,\langle N_{Tv},\psi\rangle\R).
\end{align*}                                                           
Using formula \eqref{Plancherel},  the last equality yields ${\rm Cov}\L(\L\langle\X_T(s),\varphi\R\rangle,\L\langle\X_T(t),\psi\R\rangle\R)= {\cal A} + {\cal B},$ where
			\begin{align*}
{\cal A} & =  \frac{T^3}{(2\pi)^d F^2_T}\int_{\Rd}\widehat{\varphi}(z)\overline{\widehat{\psi}(z)}
\L(\int_0^sdu\int_s^tdv+2\int_0^sdu\int_u^sdv\R)ue^{-T(v-u)|z|^{\alpha}}\,dz,\\
{\cal B} &= \frac{T^2}{(2\pi)^d F^2_T}\int_{\Rd}\widehat{\varphi}(z)\overline{\widehat{\psi}(z)}
\L(\int_0^sdu\int_s^tdv+2\int_0^sdu\int_u^sdv\R)\int_{0}^{Tu}\int_{0}^{Tu-r} e^{-[T(v+u)-2r-2w]|z|^{\alpha}}d\,{\cal U}(w)\,dr\,dz.
		\end{align*}
Notice that the term ${\cal A}$ vanishes asymptotically as $T \to \infty$
due to the choice of the norming $F_T$ and the condition $d > \alpha$.
For the second term ${\cal B}$, we perform the change of variables $r = Tr'$ and $w = Tw'$.
		\begin{align*}
			{\cal B} & =\frac{T^3}{(2\pi)^d F^2_T}\int_{\Rd}\widehat{\varphi}(z)\overline{\widehat{\psi}(z)}
			\L(\int_0^sdu\int_s^tdv+2\int_0^sdu\int_u^sdv\R)                                       			\L[
			\int_{0}^{u}\int_{0}^{u-r} e^{-T[v+u-2r-2w]|z|^{\alpha}}d\,{\cal U}(Tw)\,dr
			\R]\,dz                                                                                  \\ &=: {\cal B}_1 +{\cal B}_ 2.
		\end{align*}
After puting $y=z(T[v+u-2r-2w])^{1/{\alpha}}$, we get
		\begin{align*}{\cal B}_1 = &
			\frac{T^3{\cal U}(T)}{(2\pi)^d F^2_T T^{\frac{d}{\alpha}}}\int_{0}^s\int_{s}^t\int_{0}^u\int_{0}^{u-r}\int_{\Rd}\widehat{\varphi}\L(\frac{y}{(T[v+u-2r-2w])^{1/\alpha}}\R) 
			\overline{\widehat{\psi}\L(\frac{y}{(T[v+u-2r-2w])^{1/\alpha}}\R)}
			\\ &\cdot
			e^{-|y|^{\alpha}}(v+u-2r-2w)^{-\frac{d}{\alpha}}d\L[\frac{{\cal U}(Tw)}{{\cal U}(T)}\R]dy\,dr\,dv\,du.
		\end{align*}
Now we use that
\begin{equation}\label{Renewal_limit}
  d\L[\frac{{\cal U}(Tw)}{{\cal U}(T)}\R]\to \gamma w^{\gamma-1}\,dw\quad \mbox{and}\quad {\cal U}(T)\sim
		\frac{T^{\gamma}}{\Gamma(1+\gamma)}\quad \mbox{as}\quad  T\to\infty,
		\end{equation}
		change the order of integration (from $0\le r\le u,\ 0\le w\le u-r$ \ to \ $0\le w\le u,\ r\le u-w$) and set $F_T^2 = T^{3+\gamma-\frac{d}{\alpha}}$. Letting $T\to\infty$,
		\begin{align*}
			{\cal B}_1&\to \frac{\gamma \widehat{\varphi}(0)\overline{\widehat{\psi}}(0)\int_{\Rd}e^{-|y|^{\alpha}}\,dy}{(2\pi)
				^d\Gamma(1+\gamma)}			\int_{0}^s\int_{s}^t\int_{0}^u\int_{0}^{u-w}w^{\gamma-1}(v+u-2r-2w)^{-\frac{d}{\alpha}}dr\,dw\,dv\,du\\
						& =\frac{\gamma \widehat{\varphi}(0)\overline{\widehat{\psi}}(0)\int_{\Rd}e^{-|y|^{\alpha}}\,dy}
{2(1-\frac{d}{\alpha})(2\pi)^d\Gamma(1+\gamma)}
 \int_{0}^s\int_{s}^t\int_{0}^u
w^{\gamma-1}
			\L[
			(v+u-2w)^{1-\frac{d}{\alpha}}-(v-u)^{1-\frac{d}{\alpha}}
			\R]\,dw\,dv\,du. 
		\end{align*}
		Since
$\displaystyle
\frac{1}{(2\pi)^d}\int_{\mathbb{R}^d} e^{-|y|^{\alpha}}\,dy
= \frac{\Gamma(1/\alpha)}{\pi\alpha},
$
$
\widehat{\varphi}(0)=\langle \lambda,\varphi\rangle$ and
$
\widehat{\psi}(0)=\langle \lambda,\psi\rangle,
$
it follows that
		\begin{align*}                                                                                                                           
			{\cal B}_1 & \to \frac{\gamma\Gamma(1/\alpha)\langle\lambda,\varphi\rangle\langle\lambda,\psi\rangle}{2(1-\frac{d}{\alpha})\pi\alpha\Gamma(1+\gamma)}
			\int_{0}^s\int_{s}^t\int_{0}^u
w^{\gamma-1}
			\L[
			(v+u-2w)^{1-\frac{d}{\alpha}}-(v-u)^{1-\frac{d}{\alpha}}
			\R]\,dw\,dv\,du                                                                                                                                 =: {\cal B}_{11} +{\cal B}_{12} .
		\end{align*}
		Let $\displaystyle C_1:=\frac{\gamma\Gamma(1/\alpha)\langle\lambda,\varphi\rangle\langle\lambda,\psi\rangle}{2(1-\frac{d}{\alpha})\pi\alpha\Gamma(1+\gamma)}.$
		Using that $0\le s\le v\le t$, \  $0\le w\le u\le s$ and Fubini's theorem we obtain
		\begin{align*}
			{\cal B}_{11} & =C_1\int_{s}^t\int_{0}^s\int_{w}^sw^{\gamma-1}(v+u-2w)^{1-\frac{d}{\alpha}}\,du\,dw\,dv \\
			 &
			=\black{\frac{C_1}{2-\frac{d}{\alpha}}\int_{0}^sw^{\gamma-1}\int_{s}^t\L[
			(v+s-2w)^{2-\frac{d}{\alpha}}-(v-w)^{2-\frac{d}{\alpha}}
			\R]dv\,dw     }                                                                                              \\
			& =\frac{C_1}{(3-\frac{d}{\alpha})(2-\frac{d}{\alpha})}\int_0^sw^{\gamma-1}\L[
			(t+s-2w)^{3-\frac{d}{\alpha}}
			- (2s-2w)^{3-\frac{d}{\alpha}} - (t-w)^{3-\frac{d}{\alpha}} + (s-w)^{3-\frac{d}{\alpha}}
			\R]\,dw                                                                                                                      \\
& =			\frac{C_1}{(3-\frac{d}{\alpha})(2-\frac{d}{\alpha})}\int_0^sw^{\gamma-1}\L[
			(t+s-2w)^{3-\frac{d}{\alpha}}
			- (2^{3-\frac{d}{\alpha}}-1)(s-w)^{3-\frac{d}{\alpha}} - (t-w)^{3-\frac{d}{\alpha}}
			\R]dw,\\
			{\cal B}_{12} & =-C_1\int_{0}^s\int_{s}^t\int_{0}^uw^{\gamma-1}(v-u)^{1-\frac{d}{\alpha}}
			dw\,dv\,du                                                                                       =-\frac{C_1}{\gamma}\int_{0}^s\int_{s}^tu^{\gamma}(v-u)^{1-\frac{d}{\alpha}}dv\,du
\\
& =\black{-\frac{C_1}{\gamma(2-\frac{d}{\alpha})}\int_{0}^sw^{\gamma}\L[
			(t-w)^{2-\frac{d}{\alpha}} - (s-w)^{2-\frac{d}{\alpha}}
			\R]\,dw.}                                                                                                            
		\end{align*}
		The assumption $\alpha<d<(1+\gamma)\alpha$ implies that the above integrals are finite. We now turn to ${\cal B}_2$. Proceeding as for ${\cal B}_1,$
		\begin{align*}
			{\cal B}_2 =&\frac{2T^3{\cal U}(T)}{(2\pi)^d F^2_T T^{\frac{d}{\alpha}}}\int_0^s du\int_{u}^s dv\int_{0}^udr\int_{0}^{u-r}(v+u-2r-2w)^{-\frac{d}{\alpha}}d\L[\frac{{\cal U}(Tw)}{{\cal U}(T)}\R] \\
			& \cdot\int_{\Rd}\widehat{\varphi}\L(
			\frac{y}{T^{\frac{1}{\alpha}}(v+u-2r-2w)^{\frac{1}{\alpha}}}
			\R)
			\overline{
				\widehat{\psi}
			}\L(
			\frac{y}{T^{\frac{1}{\alpha}}(v+u-2r-2w)^{\frac{1}{\alpha}}}
			\R)e^{-|y|^{\alpha}}dy                                                                                                                                                                                  \\
			\to &\frac{2\gamma}{(2\pi)^d\Gamma(1+\gamma)}\int_{\Rd}\widehat{\varphi}(0)\overline{\widehat{\psi}}(0)e^{-|y|^{\alpha}}dy
			\int_0^s du\int_{u}^s dv\int_{0}^udr\int_{0}^{u-r}w^{\gamma-1}(v+u-2r-2w)^{-\frac{d}{\alpha}}\,dw                                                                                                   \\
			& =
			\frac{2\gamma\Gamma(1/\alpha)\langle\lambda,\varphi\rangle\langle\lambda,\psi\rangle}{\pi\alpha\Gamma(1+\gamma)}
			\int_0^s du\int_{u}^s dv\int_{0}^udw\int_{0}^{u-w}w^{\gamma-1}(v+u-2r-2w)^{-\frac{d}{\alpha}}\,dr
			\\
			& =
			2C_1
			\int_0^s du\int_{u}^s dv\int_{0}^u w^{\gamma-1}\L[
			(v+u-2w)^{1-\frac{d}{\alpha}}-(v-u)^{1-\frac{d}{\alpha}}\R]\,dw                                                                                                                             \\  =:& \ {\cal B}_{21} +{\cal B}_{22},
		\end{align*}
		where 
		\begin{align*}
			{\cal B}_{21} & =2C_1\int_0^s du\int_{u}^s dv\int_{0}^u w^{\gamma-1}
			(v+u-2w)^{1-\frac{d}{\alpha}}\,dw \\ &
			=
			2C_1\int_{0}^s dw\int_{w}^s du\int_{u}^s w^{\gamma-1}
			(v+u-2w)^{1-\frac{d}{\alpha}}\,dv                                                                                                                \\
			& = \black{\frac{2C_1}{\L(2-\frac{d}{\alpha}\R)}\int_{0}^s w^{\gamma-1}\int_{w}^s\L[
			(s+u-2w)^{2-\frac{d}{\alpha}}
			- (2u-2w)^{2-\frac{d}{\alpha}}
			\R]
			\,du\,dw}
			\\
						& =
			\frac{2C_1}{\L(2-\frac{d}{\alpha}\R)\L(3-\frac{d}{\alpha}\R)}\int_{0}^s w^{\gamma-1}
			\L[
			(2s-2w)^{3-\frac{d}{\alpha}} - (s-w)^{3-\frac{d}{\alpha}}
			-
			\frac{(2s-2w)^{3-\frac{d}{\alpha}}  }{2}
			\R] \,dw \\
			& =
						\frac{2(2^{2-\frac{d}{\alpha}}-1)C_1}{\L(2-\frac{d}{\alpha}\R)\L(3-\frac{d}{\alpha}\R)}\int_0^s
			w^{\gamma-1}(s-w)^{3-\frac{d}{\alpha}}\,dw,
\end{align*}
and
\begin{align*}{\cal B}_{22} & =-2C_1\int_0^s du\int_{u}^s dv\int_{0}^u w^{\gamma-1} (v-u)^{1-\frac{d}{\alpha}}\,dw \\ &
			=-2C_1\int_{0}^s dw\int_{w}^s du\int_{u}^sw^{\gamma-1}(v-u)^{1-\frac{d}{\alpha}}\,dv \\
			& =\black{-\frac{2C_1}{\L(2-\frac{d}{\alpha}\R)}
			\int_{0}^s w^{\gamma-1}\int_{w}^s
			{(s-u)^{2-\frac{d}{\alpha}}}
			\,du\,dw  }  \\ &
			=
			-\frac{2C_1}{\L(2-\frac{d}{\alpha}\R)\L(
				3-\frac{d}{\alpha}
				\R)}
			\int_0^s w^{\gamma-1} (s-w)^{3-\frac{d}{\alpha}}\,dw.
		\end{align*}
		In this way, we get that ${\cal B}_{11}+{\cal B}_{12}+{\cal B}_{21}+{\cal B}_{22}$ is the claimed limit covariance.

\subsection{Weak convergence via the space-time random field method}

In this section, we verify the two conditions in Theorem \ref{RFA}, which ensure the weak convergence of the sequence of processes $\{\X_T,\ T\ge1\}$ in $C([0,\Upsilon],{\cal S}'(\Rd))$ for any $\Upsilon>0$. Without loss of generality, we henceforth take $\Upsilon=1$.

We first establish the tightness of this sequence. By a classical result of Mitoma \cite[Theorem 3.1]{Mitoma}, it is sufficient to prove the tightness of the sequence of real-valued processes $\{\langle \X_T(t),\psi\rangle,\ t\ge0\}$, $T\ge1$, for each $\psi\in\E$. This will be done using the standard criterion of Billingsley \cite[Theorem 13.5]{Billingsley}.

	\subsubsection{Tightness}

	Proceeding as in the previous section,
	$$
	\mathbb{E}\L[
	\langle {\X}_T(t),\psi \rangle
	\Ph-\langle {\X}_T(s),\psi \rangle\R]^2 = {\cal D} +{\cal E},
	$$
	where
\begin{equation}\label{OtraA}
{\cal D}  =
		\frac{2T^2}{(2\pi)^dF_T^2}\int_0^s
		\int_s^ t
		\int_{\Rd}
		\L|\widehat{\psi}(z)\R|^ 2e^{-T(v-u)|z|^ {\alpha}}Tu
		\,dz\,du\,dv
\end{equation} and
$$
{\cal E}  = \frac{\gamma}{\Gamma(1+\gamma)}\int_{\Rd}e^{-|y|^{\alpha}}
		\int_{s}^t
		\int_{s}^v
		\int_{0}^u
		\int_{0}^{u-r}
		\L|\widehat{\psi}(y)\R|^2
		(v+u-2r-2w) ^{-\frac{d}{\alpha}}\,d\L[\frac{{\cal U}(Tw)}{{\cal U}(T)}\R]
		\,dr\,du\,dv\,dy .
$$
Moreover, there exist constants $c>0,$ $C>0$ such that
	\begin{align*}
		{\cal E} 
		& \le
		c
		\int_{\Rd} e^{-|y|^{\alpha}}
		\L|\widehat{\psi}(0)\R|^2dy
		\int_{s}^t
		\int_{s}^v
		\int_{0}^u
		\int_{0}^{u-r}
		(v+u-2r-2w) ^{-\frac{d}{\alpha}}
		w^{\gamma -1}\,
		dw\,dr\,du\,dv
		\\
				& =
		\frac{C}{2(1-\frac{d}{\alpha})}
		\int_s^t
		\int_s^v
		\int_0^u
		\L[
		(v+u-2w)^{1-\frac{d}{\alpha}} - (v-u)^{1-\frac{d}{\alpha}}
		\R]w^{\gamma -1}
		\,dw\,du\,dv
		\\
		& \le
		\frac{C}{2 |1-\frac{d}{\alpha} |}
		\int_s^t
		\int_s^v
		\int_0^u
		(v+u-2w)^{1-\frac{d}{\alpha}}w^{\gamma -1} \,dw\,du\,dv                                                           \\
		& =
		\frac{C}{2^{\gamma-1}4 |1-\frac{d}{\alpha} |}
		\int_s^t
		\int_s^v
		\int_0^{2u}
		(v+u-w)^{1-\frac{d}{\alpha}}w^{\gamma -1}\,dw\,du\,dv                                                             \\
		& =
		\frac{C}{2^{\gamma-1}4 |1-\frac{d}{\alpha} |}
		\int_s^t
		\int_s^v
		\int_0^{2u/(v+u)}
		(v+u)^{1-\frac{d}{\alpha}}(v+u)^{\gamma-1}(v+u)(1-w)^{1-\frac{d}{\alpha}}w^{\gamma -1}\,dw\, du\,dv.   
\end{align*}
In the following, the Beta function is denoted by {\sf B}. We proceed by
\begin{align*}
	{\cal E}	&
		\le
		\frac{C}{2^{\gamma-1}4 |1-\frac{d}{\alpha} |}
		\int_s^t
		\int_s^v
		(u+v)^{1-\frac{d}{\alpha}+\gamma}\,
		du\,dv
		\L(
		\int_0^1(1-w)^{1-\frac{d}{\alpha}}w^{\gamma-1}\,dw
		\R) \\ 
		& =
		\frac{C{\sf B}(2-\frac{d}{\alpha},\gamma)}{2^{\gamma-1}4 |1-\frac{d}{\alpha} |(2+\gamma-\frac{d}{\alpha})}
		\int_s^t\L[(2v)^{2+\gamma-\frac{d}{\alpha}}
		-
		(s+v)^{2+\gamma-\frac{d}{\alpha}}
		\R]\,
		dv                                                                                                                \\
		& =
		\frac{C{\sf B}(2-\frac{d}{\alpha},\gamma) (2+\gamma-\frac{d}{\alpha})}{2^{\gamma-1}4 |1-\frac{d}{\alpha} |(2+\gamma-\frac{d}{\alpha})}
		\int_0^t (v-s)(2\theta)^{1+\gamma-\frac{d}{\alpha}}dv,
	\end{align*}
	where the last equality follows by applying the mean value theorem to \(f(x)=x^{2+\gamma-\frac{d}{\alpha}}\), and $s+v<\theta<2v<2$. We continue by
	\begin{align}\label{Parte1} {\cal E}  & \le
		\frac{C{\sf B}(2-\frac{d}{\alpha},\gamma) }{2^{\gamma-1}4 |1-\frac{d}{\alpha} |}
		4^{1+\gamma-\frac{d}{2}}\frac{(t-s)^2}{2}
		\le
		\frac{C{\sf B}(2-\frac{d}{\alpha},\gamma) }{2^{\gamma}4 |1-\frac{d}{\alpha} |}
		4^{1+\gamma-\frac{d}{2}}
		(t-s)^{2+\gamma-\frac{d}{\alpha}}
		=
		\frac{C{\sf B}(2-\frac{d}{\alpha},\gamma) }{2^{\gamma}4 |1-\frac{d}{\alpha} |}
		4^{1+\gamma-\frac{d}{2}}
		(t-s)^{h},
	\end{align}
	where $h=2+\gamma-\frac{d}{\alpha}>1$;
	here we notice that $t-s<1$ and $\frac{d}{\alpha}-\gamma>0$, and $h>1$  because $1+\gamma -\frac{d}{\alpha}>0$.
	
	We now consider the other term, \( {\cal D} \). From \eqref{OtraA}, we obtain, by integrating with respect to \(v\) and using that \(u \le 1\),
	\begin{align*}
		{\cal D} 
		= &
		\frac{2T^2}{(2\pi)^dF_T^2}
		\int_s^ t
		\int_{\Rd}
		\L|\widehat{\psi}(z)\R|^ 2\frac{ue^{Tu|z|^ {\alpha}}}{|z|^ {\alpha}}
		\L(
		e^ {-Tu|z|^ {\alpha}} - e^ {-Tt|z|^ {\alpha}}
		\R)\,
		dz\,du                                                             \\
		& <
		\frac{2T^2}{(2\pi)^dF_T^2}
		\int_s^ t
		\int_{\Rd}
		\frac{\L|\widehat{\psi}(z)\R|^ 2}{|z|^ {\alpha}}
		\L(
		1 - e^ {-T(t-u)|z|^ {\alpha}}
		\R)\,
		dz\,du \\
		& <
		\frac{2T^2}{(2\pi)^dF_T^2}
		\int_{\Rd}
		\frac{\L|\widehat{\psi}(z)\R|^ 2}{|z|^ {\alpha}}
		\int_s^ t
		\L(
		T(t-u)|z|^ {\alpha}
		\R)^ {\delta}
		\,dz\,
		du,
\end{align*}
where $0<\delta=\gamma+1-\frac{d}{\alpha}<1.$ Hence, due to $t-u<t-s$
\begin{align}\nonumber
	{\cal D}	 < &
		\frac{2T^{2+\delta}}{(2\pi)^dT^{3+\gamma-\frac{d}{\alpha}}}
		\int_{\Rd}
		\frac{\L|\widehat{\psi}(z)\R|^ 2}{|z|^ {\alpha(1-\delta)}}dz
		\int_s^ t
		\L(
		t-s
		\R)^ {\delta}(t-s)    \\ \nonumber
		& =
		\frac{2}{(2\pi)^ d}\L(\int_{\Rd}
		\frac{|\widehat{\psi}(z)|^ 2}{|z|^ {\alpha(1-\delta)}}\,dz\R)(t-s)^ { \gamma+2-\frac{d}{\alpha}} \\ \label{Parte2} &
		=
		\frac{2}{(2\pi)^ d}\L(\int_{\Rd}
		\frac{|\widehat{\psi}(z)|^ 2}{|z|^ {\alpha(1-\delta)}}\,dz\R)(t-s)^ {h}
		.
	\end{align}
Combining the estimates \eqref{Parte1} and \eqref{Parte2}, we obtain the desired tightness property.
	\subsubsection{Weak convergence of $\{\widetilde{\X}_T\}$ in ${\cal S}'(\mathbb{R}^{d+1})$}
In this section, we verify the second condition of Theorem \ref{RFA}. To this end, we prove that the Laplace functional of $\widetilde{\X}_T$ converges, as $T \to \infty$, to the Laplace functional of a generalized process $\widetilde{\X}$. 
Finally, in the next section,  we establish \eqref{FormulaFinal}, which characterizes the limiting process $\widetilde{\X}$.
	
	\begin{theorem}\label{THM1}
		$$
		\lim_{T\to\infty}
		\EE\L[
		e^{-\L\langle\widetilde{\Phi},\widetilde{\J}_T\R\rangle
		}
		\R]
		=
		\exp
		\L\{
		\mathfrak{K}
		\int_{0}^1
		\int_{0}^{1-r}
		\int_{0}^{1-r}{\cal Q}(w,z)\varphi_2(w)\varphi_2(z)\,dz\,dw\,dl
		\R\},
		$$
		where
		$$
		\mathfrak{K}:=
		\frac
		{\gamma\L[
			\int_{\Rd}\varphi_1(x)\,dx
			\R]^2\int_{\Rd}e^{-|z|^{\alpha}}\,dz}
		{(2\pi)^d\Gamma(1+\gamma)\L(2-\frac{d}{\alpha}\R)}
		$$
		and
		\begin{equation}\label{FormulaFinal}
		\mathcal{Q}(w,z)=
		\L(\frac{d}{\alpha}-1\R)^{-1}
		\int_0^{z\wedge w}
		l^{\gamma-1}
		\L[
		(w-l)^{2-\frac{d}{\alpha}}+
		(z-l)^{2-\frac{d}{\alpha}}
		-\Ph\!
		(w+z-2l)^{2-\frac{d}{\alpha}}
		\R]\,dl.
		\end{equation}		
	\end{theorem}
	
	{\noindent\bf Proof: } 	
	Let $\wphi\in\ese.$ From \eqref{DefDelFuncional} with $\Upsilon=1$,
	\begin{align*}
		\L\langle \wphi,\widetilde{{\J}}_T\R\rangle & =
		\frac{1}{F_T}\L[
		\int_0^1\L\langle\wphi(\cdot,s),\int_0^{Ts}N_r\,dr\R\rangle\,ds
		- \int_0^1\EE\L\langle\wphi(\cdot,s),\int_0^{Ts}N_r\,dr\R\rangle\,ds
		\R]                                                                 \\
		& =
		\frac{T}{F_T}
		\L[
		\int_0^1\L\langle\wphi(\cdot,s),\int_0^{s}N_{Tr}\,dr\R\rangle\,ds
		- \int_0^1\EE\L\langle\wphi(\cdot,s),\int_0^{s}N_{Tr}\,dr\R\rangle\,ds
		\R],\quad 
		0\le r\le s\le 1.
	\end{align*}
	Using Fubini's theorem and putting
	$\Psi(x,r) = \int_r^1\wphi(x,s)\,ds,$
	$\Psi_T(x,r) = \frac{1}{F_T}\Psi(x,\frac{r}{T})$,
	\begin{align*}
	\L\langle \wphi,\widetilde{{\J}}_T\R\rangle	& =
		\frac{T}{F_T}
		\L[
		\int_{0}^1
		\int_{r}^1
		\L\langle\wphi(\cdot,s),N_{Tr}\R\rangle\,ds\,dr
		- \int_0^1
		\int_r^1
		\EE\L\langle\wphi(\cdot,s),N_{Tr}\R\rangle\,dr\,ds
		\R]                                         \\
		& =
		\frac{T}{F_T}\L[
		\int_0^1\L\langle\Psi(\cdot,r),N_{Tr}\R\rangle\,dr
		-\EE\int_0^1\L\langle\Psi(\cdot,r),N_{Tr}\R\rangle\,dr
		\R] 
\\
		& =
		\int_0^T\L\langle\Psi_T\L(\cdot,t\R),N_{t}\R\rangle\,dt - \EE\int_0^T\L\langle\Psi_T\L(\cdot,t\R),N_{t}\R\rangle\,dt.
	\end{align*}
Therefore
$$
\EE e^{-\L\langle \wphi,\widetilde{{\X}}_T\R\rangle}
=
\exp\L\{\int_0^T\EE\L\langle \Psi_T(\cdot,t),N_{t}\R\rangle\,dt\R\}
\exp
		\L\{
		\L(
		\int_{\Rd}\int_{0}^T\EE_x
		\L[
		e^{-\int_r^T\L\langle
			\Psi_T(\cdot,t),N_{t-r}(x)
			\R\rangle dt}
		\R]-1
		\R)dx \,dr
		\R\},
$$
where
$$
\int_0^T\EE\L\langle \Psi_T(\cdot,t),N_{t}\R\rangle\,dt =
\int_{\Rd}\int_0^T\int_r^T {\cal T}_{s-r}
		\Psi_T(\cdot,s)(x)\,ds\,dr\,dx.
$$
	Using that
	$$
	\int_r^T\L\langle\Psi_T(\cdot,t),N_{t-r}(x)\R\rangle\,dt = \int_0^{T-r}\L\langle\Psi_T(\cdot,t+r),N_{t}(x)\R\rangle\,dt,
	$$
	we get
	$$
	\EE_x
	\L[
	e^{-\int_r^T\L\langle
		\Psi_T(\cdot,t),N_{t-r}(x)
		\R\rangle dt}
	\R]-1 = -v_{\Psi_T}(x,r,T-r),$$ 
where 
$$v_{\Psi}(x,r,t)=\EE_x\L[
	1- e^{-\int_0^t\L\langle\Psi(\cdot,
		s+r),N_s\R\rangle ds}
	\R]
	$$ for any $\Psi\in\ese. $ Hence
	\begin{align*}
		\EE e^{-\L\langle \wphi,\widetilde{{\X}}_T\R\rangle}
		& =
		\exp\L\{
		\int_{\Rd}\int_{0}^T\int_{r}^T{\cal T}_{s-r}\Psi_T(\cdot,s)(x)\,ds\,dr\,dx
		- \int_{\Rd}\int_{0}^T
		v_{\Psi_T}(x,r,T-r)\,dr\,dx
		\R\}. 
	\end{align*}
We set $g(x)=\frac{1}{2}x^2$ and define
\[
f(x,r,t)=\EE_x\!\left[
\exp\!\left\{
-\int_0^t \Psi_T\!\left(W^{\alpha}_u,r+u\right)\,du
\right\}
\right].
\]
Recall that \(F_T^2 = T^{3+\gamma-\frac{d}{\alpha}}\) and that \(1 < \frac{d}{\alpha} < 1+\gamma\).
Applying \cite[Lemma 3.3]{ALEA}, we obtain
\begin{align*}
\int_{\mathbb{R}^d}\int_{0}^T\int_{r}^T
{\cal T}_{s-r}\Psi_T(\cdot,s)(x)\,ds\,dr\,dx
-
\int_{\mathbb{R}^d}\int_{0}^T
v_{\Psi_T}(x,r,T-r)\,dr\,dx
&= I_1(T) + I_2(T) + I_3(T) + I_4(T),
\end{align*}
where
\begin{align*}
I_1(T)
&=
\int_{\mathbb{R}^d}\int_{0}^T\int_{r}^T
{\cal T}_{s-r}\Psi_T(\cdot,s)(x)\,ds\,dr\,dx \\
&\quad -
\int_{\mathbb{R}^d}\int_{0}^T\int_{0}^{T-r}
{\cal T}_u
\bigl[
\Psi_T(\cdot,r+u)\,f(\cdot,r+u,T-r-u)
\bigr](x)\,du\,dr\,dx,
\end{align*}
\begin{align*}
-I_2(T)
&=
\frac{1}{2}
\int_{\mathbb{R}^d}\int_{0}^T
\biggl[
\int_{0}^{T-r}
{\cal T}_u
v^2_{\Psi_T}(\cdot,r+u,T-r-u)(x)\,d{\cal U}(u)
\\
&\qquad\qquad +
\int_{0}^{T-r}
\left(
\int_{0}^s {\cal T}_u \Psi_T(\cdot,T-r+u-s)(x)\,du
\right)^2
d{\cal U}(T-r-s)
\biggr] dr\,dx, \\[0.5em]
I_3(T)
&=
\frac{1}{2}
\int_{\mathbb{R}^d}\int_{0}^T\int_{0}^{T-r}
\left[
\int_{0}^s {\cal T}_u \Psi_T(\cdot,T-r+u-s)(x)\,du
\right]^2
d{\cal U}(T-r-s)\,dr\,dx, \\[0.5em]
I_4(T)
&=
\frac{1}{2}
\int_{\mathbb{R}^d}\int_{0}^T\int_{0}^{T-r}\int_{0}^{T-r-z}
{\cal T}_z
\Bigl[
\Psi_T(\cdot,r+z)\,
\mathbb{E}_{\cdot}\Bigl(
e^{-\int_{0}^u \Psi_T(W_s^{\alpha},r+z+s)\,ds}
\\
&\qquad\qquad\qquad \cdot
v^2_{\Psi_T}\bigl(
W_u^{\alpha},r+u+z,T-r-u-z
\bigr)
\Bigr)
\Bigr](x)\,
d{\cal U}(u+z)\,dz\,dr\,dx.
\end{align*}

Let \(\widetilde{\Phi}\) be of the form \(\widetilde{\Phi}(x,t)
= \varphi_1(x)\varphi_2(t)\), where \(\varphi_1\) and \(\varphi_2\) are nonnegative functions in \(\E\) and \(\mathcal{S}(\mathbb{R})\), respectively.	
	\paragraph{1. The limit $\lim_{T\to\infty} I_3(T) $}
For \(I_3(T)\), making the change of variables \(s' = T - r - s\), we obtain
\begin{align}
I_3(T)
&=
\frac{1}{2}\int_{\mathbb{R}^d}\int_{0}^T\int_{0}^{T-r}
\biggl[
\int_{0}^{T-r-s}
{\cal T}_u\bigl(\Psi_T(\cdot,s+u)\bigr)(x)\,du
\biggr]^2
\,d{\cal U}(s)\,dr\,dx \notag\\
&=
\frac{1}{2}\int_{\mathbb{R}^d}\int_{0}^T\int_{0}^{T-r}
\biggl[
\int_{s}^{T-r}
{\cal T}_{u-s}\Psi_T(\cdot,u)(x)\,du
\biggr]^2
\,d{\cal U}(s)\,dr\,dx \notag\\
&=
\frac{1}{2}\int_{\mathbb{R}^d}\int_{0}^T\int_{0}^{T-r}
\biggl[
\int_{s}^{T-r}\int_{s}^{T-r}
{\cal T}_{u-s}\Psi_T(\cdot,u)(x)\,
{\cal T}_{v-s}\Psi_T(\cdot,v)(x)\,du\,dv
\biggr]
\,d{\cal U}(s)\,dr\,dx \nonumber\\
&=
\frac{1}{2}\int_{0}^T\int_{0}^{T-r}
\biggl[
\int_{s}^{T-r}\int_{s}^{T-r}\int_{\mathbb{R}^d}
{\cal T}_{u-s}\Psi_T(\cdot,u)(x)\,
{\cal T}_{v-s}\Psi_T(\cdot,v)(x)\,dx\,du\,dv
\biggr]
\,d{\cal U}(s)\,dr. \label{edno}
\end{align}
By the Parseval–Plancherel identity,
\begin{align*}
\int_{\mathbb{R}^d}
({\cal T}_{u-s}\varphi_1)(x)\,
({\cal T}_{v-s}\varphi_1)(x)\,dx
&=
\frac{1}{(2\pi)^d}
\int_{\mathbb{R}^d}
e^{-(u+v-2s)|x|^{\alpha}}
\bigl|\widehat{\varphi}_1(x)\bigr|^2\,dx.
\end{align*}
On the other hand,
	$
	\Psi_T(\cdot,u)(x)=\frac{1}{F_T}\Psi(x,\frac{u}{T})
	=\frac{1}{F_T}\int\limits_{u/T}^1\wphi(x,s)\,ds
	=\frac{1}{F_T}\int\limits_{u/T}^1\varphi_1(x)\varphi_2(s)\,ds.
	$
	Let us write $\chi(t)=\int_t^1\varphi_2(s)\,ds$ and
	$\chi_T(t)=\chi\L(t/T\R)$. It follows from \eqref{edno} that
\begin{align*}
I_{3}(T)
&= \frac{1}{2}\int_{0}^T \int_{0}^{T-r}
\biggl[
\int_{s}^{T-r}\int_{s}^{T-r}\int_{\mathbb{R}^d}
{\cal T}_{u-s}\!\left(\frac{1}{F_T}\varphi_1\right)(x)\chi_T(u)\,
{\cal T}_{v-s}\!\left(\frac{1}{F_T}\varphi_1\right)(x)\chi_T(v)\,du\,dv\,dx
\biggr]\,d{\cal U}(s)\,dr \\
&= \frac{1}{2F_T^2(2\pi)^d}\int_{0}^T \int_{0}^{T-r}
\biggl[
\int_{s}^{T-r}\int_{s}^{T-r}\int_{\mathbb{R}^d}
e^{-(u+v-2s)|x|^{\alpha}} |\widehat{\varphi}_1(x)|^2
\chi_T(u)\chi_T(v)\,dx\,du\,dv
\biggr]\,d{\cal U}(s)\,dr.
\end{align*}
Letting \( s = Tl \), we obtain
\begin{align*}
I_{3}(T)
&= \frac{1}{2F_T^2(2\pi)^d}\int_{0}^T \int_{0}^{1-\frac{r}{T}}
\biggl[
\int_{Tl}^{T-r}\int_{Tl}^{T-r}\int_{\mathbb{R}^d}
e^{-(u+v-2Tl)|x|^{\alpha}} |\widehat{\varphi}_1(x)|^2
\chi_T(u)\chi_T(v)\,dx\,du\,dv
\biggr]\,d{\cal U}(Tl)\,dr.
\end{align*}
Applying the change of variables \( u = Tu', \; v = Tv' \), we get
\begin{align*}
I_{3}(T)
&= \frac{T^2}{2F_T^2(2\pi)^d}\int_{0}^T \int_{0}^{1-\frac{r}{T}}
\biggl[
\int_{l}^{1-\frac{r}{T}}\int_{l}^{1-\frac{r}{T}}\int_{\mathbb{R}^d}
e^{-T(u+v-2l)|x|^{\alpha}} |\widehat{\varphi}_1(x)|^2
\chi(u)\chi(v)\,dx\,du\,dv
\biggr]\,d{\cal U}(Tl)\,dr.
\end{align*}
Finally, setting \( x = z (u+v-2l)^{-1/\alpha} T^{-1/\alpha} \), we obtain
\begin{align*}
I_{3}(T)
&= \frac{T^{3-\frac{d}{\alpha}}{\cal U}(T)}{2F_T^2(2\pi)^d}
\int_{0}^1\int_{0}^{1-r}
\int_{l}^{1-r}\int_{l}^{1-r}\int_{\mathbb{R}^d}
e^{-|z|^{\alpha}}
\left|\widehat{\varphi}_1\!\left(
z (u+v-2l)^{-\frac{1}{\alpha}} T^{-\frac{1}{\alpha}}
\right)\right|^2 \\
&\qquad \cdot (u+v-2l)^{-\frac{d}{\alpha}}
\chi(u)\chi(v)\,dz\,du\,dv
\,d\!\left[\frac{{\cal U}(Tl)}{{\cal U}(T)}\right]\,dr.
\end{align*}
Since the cases \(u \leq v\) and \(v \leq u\) yield the same limit, it suffices to consider \(u \leq v\). In this case, \(0 \leq l \leq u \leq v \leq 1 - r\). Hence,
\begin{align}\nonumber
\lim_{T\to\infty} I_3(T)
&=
2\lim_{T\to\infty}
\frac{T^{3-\frac{d}{\alpha}} \, {\cal U}(T)}{2(2\pi)^d F_T^2}
\int_{0}^1 \int_{0}^{1-r} \int_{0}^v
\biggl(
\int_{0}^u \int_{\mathbb{R}^d}
e^{-|z|^{\alpha}}
\left|
\widehat{\varphi}_1\!\left(
z (u+v-2l)^{-\frac{1}{\alpha}} T^{-\frac{1}{\alpha}}
\right)
\right|^2  \\ \nonumber
&\qquad\qquad \cdot
(u+v-2l)^{-\frac{d}{\alpha}}\,
d\!\left[\frac{{\cal U}(Tl)}{{\cal U}(T)}\right]
\,dz
\biggr)
\chi(u)\chi(v)\,du\,dv\,dr \label{Limit_of_I(3)} \\
&=
\frac{\gamma |\widehat{\varphi}_1(0)|^2}{(2\pi)^d \Gamma(1+\gamma)}
\left(
\int_{\mathbb{R}^d} e^{-|z|^{\alpha}}\,dz
\right)
\int_{0}^1 \int_{0}^{1-r} \int_{0}^v \int_{0}^u
(u+v-2l)^{-\frac{d}{\alpha}} l^{\gamma-1}
\chi(u)\chi(v)\,dl\,du\,dv\,dr,
\end{align}
where we used \eqref{Renewal_limit} to obtain the last equality.
	\paragraph{2. The limit $\lim_{T\to\infty} I_1(T) $}
We have
\begin{align*}
I_1(T)
&=
\int_{\mathbb{R}^d}\int_{0}^T\int_{r}^T
{\cal T}_{s-r}\Psi_T(\cdot,s)(x)\,ds\,dr\,dx \\
&\quad -
\int_{\mathbb{R}^d}\int_{0}^T\int_{0}^{T-r}
{\cal T}_{u}
\bigl[
\Psi_T(\cdot,r+u)f(\cdot,r+u,T-r-u)
\bigr](x)\,du\,dr\,dx.
\end{align*}

Performing the change of variables \(u = s - r\), we obtain
\begin{align*}
I_1(T)
&=
\int_{\mathbb{R}^d}\int_{0}^T\int_{0}^{T-r}
{\cal T}_u \Psi_T(\cdot,r+u)(x)\,du\,dr\,dx \\
&\quad -
\int_{\mathbb{R}^d}\int_{0}^T\int_{0}^{T-r}
{\cal T}_u
\bigl[
\Psi_T(\cdot,r+u)f(\cdot,r+u,T-r-u)
\bigr](x)\,du\,dr\,dx \\
&=
\int_{\mathbb{R}^d}\int_{0}^T\int_{0}^{T-r}
{\cal T}_u
\bigl[
\Psi_T(\cdot,r+u)\bigl(1 - f(\cdot,r+u,T-r-u)\bigr)
\bigr](x)\,du\,dr\,dx.
\end{align*}

Using the invariance of the semigroup \({\cal T}\) with respect to the Lebesgue measure, we get
\begin{align}
I_1(T)
&=
\int_{\mathbb{R}^d}\int_{0}^T\int_{0}^{T-r}
\Psi_T(x,r+u)\bigl(1 - f(x,r+u,T-r-u)\bigr)
\,du\,dr\,dx \nonumber \\
&=
\int_{\mathbb{R}^d}\int_{0}^T\int_{0}^{T-r}
\Psi_T(x,r+u)\int_{0}^{T-r-u}
{\cal T}_s
\bigl[
\Psi_T(\cdot,r+u+s)f(\cdot,r+u+s,T-r-u-s)
\bigr](x)\,ds\,du\,dr\,dx \nonumber \\
&\le
\int_{\mathbb{R}^d}\int_{0}^T\int_{r}^{T}
\Psi_T(x,u)\int_{0}^{T-u}
({\cal T}_s \Psi_T(\cdot,u+s))(x)\,ds\,du\,dr\,dx,
\label{dve}
\end{align}
Let \(\mathfrak{c} := \sup_{t \geq 0} \varphi_2(t)\). Then
\[
\Psi_T(x,s)
= \frac{1}{F_T}\varphi_1(x)\int_{s/T}^{1}\varphi_2(t)\,dt
\leq \frac{1}{F_T}\varphi_1(x)\int_{0}^{1}\varphi_2(t)\,dt
\leq \frac{\mathfrak{c}}{F_T}\varphi_1(x).
\]

From \eqref{dve} and the above bound, we obtain
\begin{align*}
I_1(T)
&\leq \frac{\mathfrak{c}^2}{F_T^2}
\int_{\mathbb{R}^d}\int_{0}^T\int_{0}^T
\varphi_1(x)
\left(
\int_{0}^{T-u} {\cal T}_s \varphi_1(x)\,ds
\right)\,du\,dr\,dx \\
&\leq \frac{\mathfrak{c}^2 T^2}{F_T^2}
\int_{\mathbb{R}^d}
\varphi_1(x)
\left(
\int_{0}^T {\cal T}_s \varphi_1(x)\,ds
\right)\,dx.
\end{align*}

Using Parseval's identity,
\begin{align*}
I_1(T)
&= \frac{\mathfrak{c}^2 T^2}{(2\pi)^d F_T^2}
\int_{\mathbb{R}^d}\int_{0}^T
|\widehat{\varphi}_1(x)|^2 e^{-s|x|^{\alpha}}\,ds\,dx 
= \frac{\mathfrak{c}^2 T^2}{(2\pi)^d F_T^2}
\int_{\mathbb{R}^d}
|\widehat{\varphi}_1(x)|^2
\frac{1 - e^{-T|x|^{\alpha}}}{|x|^{\alpha}}\,dx.
\end{align*}

Recall the elementary inequality \(1 - e^{-x} \leq x^{\delta}\) for \(x \geq 0\) and \(0 \leq \delta \leq 1\).
Choose \(0 < \delta < 1\) such that
\[
\delta < 1 + \gamma - \frac{d}{\alpha},
\qquad
\delta < \frac{d}{\alpha} - 1.
\]

Then
\begin{align*}
I_1(T)
&\leq
\frac{\mathfrak{c}^2 T^2}{(2\pi)^d F_T^2}
\int_{\mathbb{R}^d}
|\widehat{\varphi}_1(x)|^2
\frac{T^{\delta} |x|^{\alpha \delta}}{|x|^{\alpha}}\,dx 
=
\frac{\mathfrak{c}^2}{(2\pi)^d \, T^{1+\gamma-\frac{d}{\alpha}-\delta}}
\int_{\mathbb{R}^d}
\frac{|\widehat{\varphi}_1(x)|^2}{|x|^{\alpha(1-\delta)}}\,dx.
\end{align*}

Since \(\alpha(1-\delta) < d\), the last integral is finite, and as \(1+\gamma-\frac{d}{\alpha}-\delta > 0\), it follows that
\[
\lim_{T\to\infty} I_1(T) = 0.
\]
	\paragraph{3. Convergence of \(I_4(T)\).}

Clearly, \(0 \leq -I_4(T)\), and
\begin{align*}
-I_4(T)
&=
\frac{1}{2}
\int_{\mathbb{R}^d}\int_{0}^T\int_{0}^{T-r}\int_{0}^{T-r-z}
{\cal T}_z
\Bigl[
\Psi_T(\cdot,r+z)\,
\mathbb{E}_{\cdot}\Bigl(
e^{-\int_{0}^u \Psi_T(W_s^{\alpha},r+z+s)\,ds}
\, v^2_{\Psi_T}\!\left(
W_u^{\alpha}, r+u+z, T-r-u-z
\right)
\Bigr)
\Bigr](x) \\
&\qquad\qquad \cdot d{\cal U}(u+z)\,dz\,dr\,dx.
\end{align*}

Setting \(v = T - r - u - z\), we obtain
\begin{align*}
-I_4(T)
&\leq
\frac{1}{2}
\int_{\mathbb{R}^d}\int_{0}^T\int_{0}^{T-r}\int_{0}^{T-r-z}
{\cal T}_z
\Bigl[
\Psi_T(\cdot,r+z)\,
\mathbb{E}_{\cdot}\Bigl(
v^2_{\Psi_T}\bigl(
W^{\alpha}_{T-r-v-z}, T-v, v
\bigr)
\Bigr)
\Bigr](x) \\
&\qquad\qquad \cdot d{\cal U}(T-r-v)\,dz\,dr\,dx.
\end{align*}
From (4.23) in \cite{ALEA} and the definition of \(\Psi_T\), it follows that
\begin{align*}
-I_4(T)
&\leq
\frac{\mathfrak{c}}{2F_T^3}
\int_{\mathbb{R}^d}\int_{0}^T\int_{0}^{T-r}\int_{0}^{T-r-z}
{\cal T}_z
\Bigl[
\varphi_1(\cdot)
\Bigl(
\int_{0}^T {\cal T}_l \varphi_1(\cdot)\,dl
\Bigr)^2
\Bigr](x)
\,d{\cal U}(T-r-v)\,dz\,dr\,dx \\
&\leq
\frac{\mathfrak{c}}{2F_T^3}
\int_{\mathbb{R}^d}\int_{0}^T\int_{0}^{T-r}\int_{0}^{T-z}
{\cal T}_z
\Bigl[
\varphi_1(\cdot)
\Bigl(
\int_{0}^T {\cal T}_l \varphi_1(\cdot)\,dl
\Bigr)^2
\Bigr](x)
\,d{\cal U}(T-r-v)\,dz\,dr\,dx.
\end{align*}
Arguing as in \cite[p.~612]{ALEA}, we obtain, for some constant \(C = C(\alpha,\varphi_1,d)>0\),
\[
0 \leq -I_4(T)
\leq
\frac{C T^{\gamma} T}{F_T^3}
=
\frac{C T^{1+\gamma}}{T^{\frac{3}{2}(3+\gamma-\frac{d}{\alpha})}}
=
\frac{C}{
T^{\frac{7}{2}+\frac{\gamma}{2}-\frac{3}{2}\frac{d}{\alpha}}
}.
\]
Noting that
\[
\frac{7}{2}+\frac{\gamma}{2}-\frac{3}{2}\frac{d}{\alpha}
>
\frac{7}{2}+\frac{\gamma}{2}-\frac{3}{2}(1+\gamma)
=
2-\gamma > 0,
\]
we conclude that
\[
\lim_{T \to \infty} I_4(T) = 0.
\]

\paragraph{4. Convergence of $I_2(T)$}

We start from the following expression for $I_2(T)$:
\begin{align*}
0 \le -2 I_2(T)
&= -\int_{\mathbb{R}^d} \int_0^T \int_0^{T-r}
v^2_{\Psi_T}(x,s,T-r-s)\, d\mathcal{U}(s)\, dr\, dx \\
&\quad +
\int_{\mathbb{R}^d} \int_0^T \int_0^{T-r}
\left(
\int_0^s \mathcal{T}_u\big[
\Psi_T(\cdot,T-r+u-s)
\big](x)\, du
\right)^2
d\mathcal{U}(T-r-s)\, dr\, dx.
\end{align*}

Performing the change of variables $s' = T-r-s$ in the first integral, we obtain
\begin{align}
-2 I_2(T)
&= -\int_{\mathbb{R}^d} \int_0^T \int_0^{T-r}
v^2_{\Psi_T}(x,T-r-s,s)\, d\mathcal{U}(T-r-s)\, dr\, dx \notag \\
&\quad +
\int_{\mathbb{R}^d} \int_0^T \int_0^{T-r}
\left(
\int_0^s \mathcal{T}_u\big[
\Psi_T(\cdot,T-r+u-s)
\big](x)\, du
\right)^2
d\mathcal{U}(T-r-s)\, dr\, dx \notag \\
&=
\int_{\mathbb{R}^d} \int_0^T \int_0^{T-r}
\Bigg[
\left(
\int_0^s \mathcal{T}_u\big[
\Psi_T(\cdot,T-r+u-s)
\big](x)\, du
\right)^2
- v^2_{\Psi_T}(x,T-r-s,s)
\Bigg]
d\mathcal{U}(T-r-s)\, dr\, dx.
\label{tri}
\end{align}

Next, by inequality (4.29) in \cite{ALEA}, we have
\begin{align}
0 \le\,
&\int_0^s \mathcal{T}_u\big(
\Psi_T(\cdot,T-r+u-s)
\big)\, du
- v_{\Psi_T}(x,T-r-s,s) \notag \\
\le\,
&\int_0^s \mathcal{T}_u \Bigg[
\Psi_T(\cdot,T-r+u-s)
\left(
\int_0^{s-u} \mathcal{T}_w
\Psi_T(\cdot,T-r-s+u+w)\, dw
\right)
\Bigg](x)\, du
\label{chetiri} \\
&\quad +
\frac{1}{2}
\int_0^s \mathcal{T}_u \Bigg[
\int_0^{s-u} \mathcal{T}_w
\Psi_T(\cdot,T-r-s+u+w)\, dw
\Bigg]^2 (x)\, d\mathcal{U}(u). \notag
\end{align}

On the other hand, by relation (4.23) in \cite{ALEA}, we have
\[
v_{\Psi_T}(x,T-r-s,s)
\le
\int_0^s \mathcal{T}_u \big[
\Psi_T(\cdot,T-r-s+u)
\big](x)\, du.
\]

Consequently,
\begin{align*}
0 \le\,
&\left(
\int_0^s \mathcal{T}_u \Psi_T(\cdot,T+u-r-s)(x)\, du
\right)^2
- v^2_{\Psi_T}(x,T-r-s,s) \\
=\,
&\left(
\int_0^s \mathcal{T}_u \Psi_T(\cdot,T+u-r-s)(x)\, du
+ v_{\Psi_T}(x,T-r-s,s)
\right) \\
&\quad \times
\left(
\int_0^s \mathcal{T}_u \Psi_T(\cdot,T+u-r-s)(x)\, du
- v_{\Psi_T}(x,T-r-s,s)
\right).
\end{align*}

It follows that
\begin{align*}
&\left(
\int_0^s \mathcal{T}_u \Psi_T(\cdot,T+u-r-s)(x)\, du
\right)^2
- v^2_{\Psi_T}(x,T-r-s,s) \\
&\quad \le
2
\int_0^s \mathcal{T}_u \Psi_T(\cdot,T+u-r-s)(x)\, du \\
&\qquad \times
\left[
\int_0^s \mathcal{T}_u \Psi_T(\cdot,T+u-r-s)(x)\, du
- v_{\Psi_T}(x,T-r-s,s)
\right] \\
&\quad \le
2\int_0^s \mathcal{T}_u \Psi_T(\cdot,T+u-r-s)(x)\, du \\
&\qquad \times
\Bigg[
\int_0^s \mathcal{T}_u \Bigg\{
\Psi_T(\cdot,T-r+u-s)
\int_0^{s-u} \mathcal{T}_w \Psi_T(\cdot,T-r-s+u+w)\, dw
\Bigg\}(x)\, du \\
&\qquad\qquad +
\frac{1}{2}
\int_0^s \mathcal{T}_u \left(
\int_0^{s-u} \mathcal{T}_w \Psi_T(\cdot,T-r-s+u+w)\, dw
\right)^2(x)\, d\mathcal{U}(u)
\Bigg].
\end{align*}	where the  last inequality follows from \eqref{chetiri}. Hence, from \eqref{tri} we obtain
	\begin{align*}
		0  \le&-2I_2(T)                                                                 \\ \le &
		\int_{\Rd}\int_{0}^T\int_{0}^{T-r}\L[
		\L(
		\int_0^s{\cal T}_u\Psi_T(\cdot,T+u-r-s) \,du\R)\R. \\ 
&\L.\times\int_0^s{\cal T}_u\L(
		\int_0^{s-u}{\cal T}_w\Psi_T(\cdot,T-r-s+u+w)\,dw \R)^2 \,d{\cal U}(u)\R.
		\\
		& +2\L.\int_0^s{\cal T}_u\Psi_T(\cdot,T+u-r-s) \,du \R. \\ 
&\L.\times\int_0^s\L\{
		{\cal T}_u\Psi_T(\cdot,T-r+u-s)\int_0^{s-u}{\cal T}_w\Psi_T(\cdot,T-r-s+u+w)\,dw
		\R\}
		\,du
		\R](x)
		\, d{\cal U}(T-r-s)\,dr\,dx.\\
		:=&\int_{\Rd}\int_{r=0}^T\int_{s=0}^{T-r}\L[\mathfrak{I}_1(T)\mathfrak{R}_1(T)+\Ph
		2 \,\mathfrak{I}_1(T) \mathfrak{R}_2(T)
		\R]\,d{\cal U}(T-r-s)\,dr\,dx,
	\end{align*}
where
\begin{align*}
\mathfrak{I}_1(T) & = \int_0^s{\cal T}_u\Psi_T(\cdot,T+u-r-s) \,du, \\
\mathfrak{R}_1(T) & = \int_0^s{\cal T}_u\L(
		\int_0^{s-u}{\cal T}_w\Psi_T(\cdot,T-r-s+u+w)\,dw \R)^2 \,d{\cal U}(u),\\
\mathfrak{R}_2(T) & = \int_0^s\L\{
		{\cal T}_u\Psi_T(\cdot,T-r+u-s)\int_0^{s-u}{\cal T}_w\Psi_T(\cdot,T-r-s+u+w)\,dw
		\R\}
		\,du.
\end{align*}
	Let us denote
	\begin{align*}
		H(T)   & =\int_{\Rd}\int_0^T\int_0^{T-r}
		\L(
		\int_0^s{\cal T}_u\Psi_T(\cdot,T+u-r-s)(x)\,du\R)^2
		\,d{\cal U}(T-r-s)\,dr\,dx \\ & =	
		\int_{\Rd}\int_{r=0}^T\int_{s=0}^{T-r} \mathfrak{I}_1^2(T)\,d{\cal U}(T-r-s)\,dr\,dx,  \\
		R_1(T) & =
		\int_{\Rd}\int_{r=0}^T\int_{s=0}^{T-r}\L[
		\int_0^s{\cal T}_u
		\Psi_T(\cdot,T-r-s+u)\int_0^{s-u}{\cal T}_w\Psi_T(\cdot,T-r-s+u+w)\,dw
		\R]^2d{\cal U}(T-r-s)\,dr\,dx \\ 
		&\phantom{XXX}=\int_{\Rd}\int_0^T\int_0^{T-r}\mathfrak{R}^2_1(T)\,d{\cal U}(T-r-s)\,dr\,dx,  \\
		R_2(T) &=
		\,d{\cal U}(T-r-s)\,dr\,dx, 
		\L[
		\int_0^s{\cal T}_u\L(\int_0^{s-u}
		{\cal T}_w\Psi_T(\cdot,T-s-r+u+w)\,dw\R)^2\,d{\cal U}(u)
		\R]^2\,d{\cal U}(T-r-s)\,dr\,dx\\
		&\phantom{XXXX} = \int_{\Rd}\int_0^T\int_0^{T-r}\mathfrak{R}_2^2(T)
		\,d{\cal U}(T-r-s)\,dr\,dx.
	\end{align*}
	Then, applying the Cauchy-Schwarz inequality to the measure $
	\int_{\Rd}\int_{r=0}^T\int_0^{T-r}{\cal U}(T-r-s)\,dr\,dx,
	$
	we obtain
	$$
	-2I_2(T)\le \mathrm{Const.}\sqrt{H(T)}\L(\sqrt{R_1(T)} +\sqrt{R_2(T)}\R).$$
	Proceeding similarly as in \cite[(3.14)]{BGTAnnals}, one can show that above $H$ is bounded.
	We now show that both $R_1(T)$ and $R_2(T)$ vanish as $T \to \infty$. 
We begin with $R_1(T)$, for which we have
\begin{align*}
R_1(T)
&=
\int_{\mathbb{R}^d}\int_0^T\int_0^{T-r}
\left\{
\int_0^s \mathcal{T}_u \Big[
\Psi_T(\cdot,T-r+u-s)
\int_0^{s-u} \mathcal{T}_w \Psi_T(\cdot,T-r-s+u+w)\, dw
\Big](x)\, du
\right\}^2 \\
&\qquad \times d\mathcal{U}(T-r-s)\, ds\, dr\, dx.
\end{align*}

We perform the rescaling $r = Tr'$ and $s = Ts'$, which yields
\begin{align*}
R_1(T)
&\le
\frac{T^2}{F_T^4}
\int_{\mathbb{R}^d}\int_0^1\int_0^{1-r}
\left[
\int_0^{Ts} \mathcal{T}_u \left(
\varphi_1(\cdot)
\int_0^{Ts-u} \mathcal{T}_w \varphi_1(\cdot)\, dw
\right)(x)\, du
\right]^2 \\
&\qquad \times d\mathcal{U}(T-Tr-Ts)\, dr\, ds\, dx.
\end{align*}

Next, introducing the change of variables $u = Tu'$ and $w = Tw'$, we obtain
\begin{align*}
R_1(T)
&\le
\frac{T^6}{F_T^4}
\int_{\mathbb{R}^d}\int_0^1\int_0^{1-r}
\left(
\int_0^s \mathcal{T}_{Tu}
\left(
\varphi_1(\cdot)
\int_0^{s-u} \mathcal{T}_{Tw} \varphi_1(\cdot)\, dw
\right)(x)\, du
\right)^2 \\
&\qquad \times d\mathcal{U}\big(T(1-r-s)\big)\, dr\, ds\, dx.
\end{align*}

Estimating the integrals uniformly, we further obtain
\begin{align*}
R_1(T)
&\le
\frac{T^6}{F_T^4}
\int_{\mathbb{R}^d}\int_0^1\int_0^{1-r}
\left[
\int_0^1 \mathcal{T}_{Tu} \varphi_1(x)
\left(
\int_0^1 \mathcal{T}_{Tw} \varphi_1(\cdot)\, dw
\right)(x)\, du
\right]^2
d\mathcal{U}\big(T(1-r-s)\big)\, dr\, ds\, dx \\
&\le
\frac{T^6\, u(T)}{F_T^4}
\int_{\mathbb{R}^d}
\left[
\int_0^1 \mathcal{T}_{Tu} \varphi_1(x)
\left(
\int_0^1 \mathcal{T}_{Tw} \varphi_1(\cdot)\, dw
\right)(x)\, du
\right]^2 dx.
\end{align*}
Finally, using the self-similarity  of the $\alpha$-stable semigroup and the scaling properties assumptions, we conclude that
\[
R_1(T)
\le
\mathrm{Const.}\, \frac{T^{6+\gamma}}{T^{2(3+\gamma-\frac{d}{\alpha})+\frac{3d}{\alpha}}}
=
\frac{\mathrm{Const.}}{T^{\gamma + \frac{d}{\alpha}}}
\;\xrightarrow[T\to\infty]{}\; 0,
\]
which proves the desired convergence.	

	We next consider $R_2(T)$. We have
\begin{align*}
R_2(T)
\le\,
&\int_{\mathbb{R}^d}\int_0^T\int_0^{T-r}
\left[
\int_0^s \mathcal{T}_u \left(
\int_0^{s-u} \mathcal{T}_w \Psi_T(\cdot)\, dw
\right)^2 (x)\, d\mathcal{U}(u)
\right]^2 d\mathcal{U}(T-r-s)\, dr\, dx.
\end{align*}

We perform the rescaling $u=Tu'$, $r=Tr'$, and $s=Ts'$, which yields
\begin{align*}
R_2(T)
\le\,
&\frac{T}{F_T^4}
\int_{\mathbb{R}^d}\int_0^1\int_0^{1-r}
\left[
\int_0^s \mathcal{T}_{Tu}
\left(
\int_0^{Ts-Tu} \mathcal{T}_w \varphi_1(\cdot)\, dw
\right)^2 (x)\, d\mathcal{U}(Tu)
\right]^2
d\mathcal{U}(T(1-r-s))\, dr\, ds\, dx\\
\le\,
&\frac{T^5}{F_T^4}
\int_{\mathbb{R}^d}\int_0^1\int_0^{1-r}
\left[
\int_0^1 \mathcal{T}_{Tu}
\left(
\int_0^1 \mathcal{T}_{Tw} \varphi_1(\cdot)\, dw
\right)^2 (x)\, d\mathcal{U}(Tu)
\right]^2
d\mathcal{U}(T(1-r-s))\, dr\, ds\, dx.
\end{align*}

Estimating uniformly in $r$, we further obtain
\begin{align*}
R_2(T)
\le\,
\frac{T^5 {\cal U}(T)^3}{F_T^4}
\int_{\mathbb{R}^d}
\left[
\int_0^1 \mathcal{T}_{Tu}
\left(
\int_0^1 \mathcal{T}_{Tw} \varphi_1(\cdot)\, dw
\right)^2 (x)\, \frac{d\mathcal{U}(Tu)}{{\cal U}(T)}
\right]^2 dx.
\end{align*}
Using the integral representation of the semigroup, we rewrite
\begin{align*}
R_2(T)
\le\,
&\frac{T^5 {\cal U}(T)^3}{F_T^4}
\int_{\mathbb{R}^d}
\Bigg[
T^{-\frac{d}{\alpha}}
\int_{\mathbb{R}^d}\int_0^1
p_u\big((x-y)T^{-1/\alpha}\big) \\
&\qquad \times
\left(
\int_{\mathbb{R}^d}\int_0^1
T^{-\frac{d}{\alpha}} p_w\big((y-z)T^{-1/\alpha}\big)
\varphi_1(z)\, dz\, dw
\right)^2
\frac{d\mathcal{U}(Tu)}{{\cal U}(T)}\, dy
\Bigg]^2 dx.
\end{align*}
We now perform the change of variables
$x = T^{1/\alpha} x'$, $y = T^{1/\alpha} y'$, $z = T^{1/\alpha} z'$. Then
\begin{align*}
R_2(T)
=\,
\mathrm{Const.}\,
\frac{1}{T^{1-\gamma+\frac{d}{\alpha}}}
\left\|
f_{1,T} * (r * g_{1,T})^2
\right\|_2^2
\xrightarrow[T\to\infty]{} 0,
\end{align*}
where
\[
f_{1,T}(x)=\int_0^1 p_u(x)\frac{d\mathcal{U}(Tu)}{{\cal U}(T)}, \quad
g_{1,T}(x)=T^{\frac{d}{\alpha}}\varphi_1(T^{1/\alpha}x),
\quad
r(x)=\int_0^1 p_u(x)\, du,
\]
and $\|f_{1,T}\|_1 \le\mathrm{Const.}$ for $T$ sufficiently large, and we have used Young’s inequality,
\[
\|f_{1,T} * (r * g_{1,T})^2\|_2^2
\le
\|f_{1,T}\|_1^2 \, \|r\|_2^4 \, \|g_{1,T}\|_1^4
< \infty.
\]	

\paragraph{5.  Proof of \eqref{FormulaFinal}}
	From \eqref{Limit_of_I(3)} we know that
\begin{align*}
\lim_{T\to\infty}\mathbb{E}\left[
e^{-\langle \widetilde{\Phi}, \widetilde{\mathcal{J}}_T\rangle}
\right]
&=
\exp\Bigg\{
\frac{\beta\gamma\big[\int_{\mathbb{R}^d}\varphi_1(x)\,dx\big]^2}
{(2\pi)^d\Gamma(1+\gamma)}
\left(\int_{\mathbb{R}^d}e^{-|z|^{\alpha}}dz\right)
\\
&\qquad \times
\int_0^1 \int_0^{1-r} \int_0^{1-r} \int_0^{u\wedge v}
l^{\gamma-1}(u+v-2l)^{-\frac{d}{\alpha}}
\chi(u)\chi(v)\, dl\, du\, dv\, dr
\Bigg\},
\end{align*}
where $\chi(z)=\int_z^1 \varphi_2(u)\,du$, $0\le z\le 1$.

For convenience, we set
\[
M :=
\int_0^1 \int_0^{1-r} \int_l^{1-r} \int_l^{1-r}
l^{\gamma-1}(u+v-2l)^{-\frac{d}{\alpha}}
\chi(u)\chi(v)\, du\, dv\, dl\, dr.
\]

Expanding $\chi$, we rewrite $M$ as
\begin{align*}
M
&=
\int_0^1 \int_0^{1-r} \int_l^{1-r} \int_l^{1-r}
\int_l^{1-r} \int_l^{1-r}
l^{\gamma-1}(u+v-2l)^{-\frac{d}{\alpha}} \\
&\qquad \times \varphi_2(w)\varphi_2(z)\,
dw\, dz\, du\, dv\, dl\, dr.
\end{align*}

We now split the integration domain according to the order of $u$ and $v$.
Assume first that $z<w$.

\medskip
\noindent
\textbf{Case 1: $u<v$.} In this case,
\begin{align*}
M_1
&=
\int_0^1 \int_0^{1-r} \int_l^{1-r} \int_l^{1-r}
\left(
\int_l^z \int_u^w
l^{\gamma-1}(u+v-2l)^{-\frac{d}{\alpha}}\, dv\, du
\right)
\varphi_2(w)\varphi_2(z)\, dw\, dz\, dl\, dr.
\end{align*}

Computing the inner integral explicitly, we obtain
\begin{align*}
M_1
&=
\int_0^1 \int_0^{1-r} \int_l^{1-r} \int_l^{1-r}
\frac{l^{\gamma-1}}{\big(1-\frac{d}{\alpha}\big)\big(2-\frac{d}{\alpha}\big)} \\
&\qquad \times
\Big[
(z+w-2l)^{2-\frac{d}{\alpha}}
-(w-l)^{2-\frac{d}{\alpha}}
-\tfrac{1}{2}(2z-2l)^{2-\frac{d}{\alpha}}
\Big] \\
&\qquad \times \varphi_2(w)\varphi_2(z)\, dw\, dz\, dl\, dr.
\end{align*}

\medskip
\noindent
\textbf{Case 2: $u>v$.} Similarly,
\begin{align*}
M_2
&=
\int_0^1 \int_0^{1-r} \int_l^{1-r} \int_l^{1-r}
\frac{l^{\gamma-1}}{\big(1-\frac{d}{\alpha}\big)\big(2-\frac{d}{\alpha}\big)} \\
&\qquad \times
\Big[
\tfrac{1}{2}(2z-2l)^{2-\frac{d}{\alpha}}
-(z-l)^{2-\frac{d}{\alpha}}
\Big]
\varphi_2(w)\varphi_2(z)\, dw\, dz\, dl\, dr.
\end{align*}

Combining both contributions, we deduce that
\begin{align*}
M
&= M_1 + M_2 \\
&=
\frac{1}{\big(1-\frac{d}{\alpha}\big)\big(2-\frac{d}{\alpha}\big)}
\int_0^1 \int_0^{1-r} \int_l^{1-r} \int_l^{1-r}
l^{\gamma-1}
(z-l)^{2-\frac{d}{\alpha}}
\varphi_2(w)\varphi_2(z)\, dw\, dz\, dl\, dr.
\end{align*}

The case $z>w$ can be treated analogously.	
	
	\section{Proof of Theorem \ref{Properties}}\label{ProofProperties} 

Recall the covariance representation \eqref{F2}:
\begin{align*}
\mathbb{E}\big[\X(s)\X(t)\big]
=
\C \int_0^s w^{\gamma-1}
\Big(
3(s-w)^{3-\da}
+ (t-w)^{3-\da}
+ \frac{(3-\da)w}{\gamma}\big[(t-w)^{2-\da}-(s-w)^{2-\da}\big]
- (t+s-2w)^{3-\da}
\Big)\,dw,
\end{align*}
where
\[
\C=\frac{\gamma\Gamma(1/\alpha)\langle\lambda,\varphi\rangle\langle\lambda,\psi\rangle}
{2\big(\frac{d}{\alpha}-1\big)\big(2-\frac{d}{\alpha}\big)\big(3-\frac{d}{\alpha}\big)\pi\alpha\Gamma(1+\gamma)}.
\]
	
		\subsection{Proof of Theorem \ref{Properties}.1}

It follows that, for any $a>0$,
\begin{align*}\lefteqn{
\mathbb{E}\big[\X(as)\X(at)\big]}\\ &
=
\C \int_0^{as} w^{\gamma-1}
\Big(
3(as-w)^{3-\da}
+ (at-w)^{3-\da}
+ \frac{(3-\da)w}{\gamma}\big[(at-w)^{2-\da}-(as-w)^{2-\da}\big]
- (as+at-2w)^{3-\da}
\Big)\,dw.
\end{align*}
Performing the change of variables $w=aw'$, a straightforward computation shows that each term scales as
\[
\int_0^{as} w^{\gamma-1}(as-w)^{3-\da}\,dw
= a^{3+\gamma-\da}\int_0^s w^{\gamma-1}(s-w)^{3-\da}\,dw,
\]
and similarly for the remaining terms. Consequently,
\[
\mathbb{E}\big[\X(as)\X(at)\big]
= a^{3+\gamma-\da}\,\mathbb{E}\big[\X(s)\X(t)\big].
\]
This establishes the self-similarity of $\X$ with index $\tfrac{3+\gamma-\da}{2}$.

\subsection{Proof of Theorem \ref{Properties}.2}

For simplicity, set
\[
Q(s,t)=\mathbb{E}\big[\X(s)\X(t)\big].
\]
Using \eqref{F2}, differentiation with respect to $s$ yields
\begin{align*}
\frac{\partial Q(s,t)}{\partial s}
&=
\C\Bigg\{
(3-\tfrac{d}{\alpha})\gamma^{-1}
s^{\gamma}(t-s)^{2-\frac{d}{\alpha}} \\
&\quad +
\int_0^s w^{\gamma-1}
\Big[
3(3-\tfrac{d}{\alpha})(s-w)^{2-\frac{d}{\alpha}}
- (3-\tfrac{d}{\alpha})\gamma^{-1}(2-\tfrac{d}{\alpha})w(s-w)^{1-\frac{d}{\alpha}} \\
&\qquad\qquad\quad
- (3-\tfrac{d}{\alpha})(t+s-2w)^{2-\frac{d}{\alpha}}
\Big]\,dw
\Bigg\}.
\end{align*}
Rewriting, we obtain the scaling form
\[
\frac{\partial Q(s,t)}{\partial s}
=
\C\,(3-\tfrac{d}{\alpha})\, s^{2+\gamma-\frac{d}{\alpha}}
\, f\!\left(\frac{t}{s}\right),
\quad s>0,
\]
where, for $x\ge 1$,
\begin{align*}
f(x)
&=
\gamma^{-1}(x-1)^{2-\frac{d}{\alpha}} \\
&\quad +
\int_0^1 w^{\gamma-1}
\Big[
3(1-w)^{2-\frac{d}{\alpha}}
-\gamma^{-1}(2-\tfrac{d}{\alpha})w(1-w)^{1-\frac{d}{\alpha}}
-(x+1-2w)^{2-\frac{d}{\alpha}}
\Big]\,dw.
\end{align*}
In particular, setting $\displaystyle g(s,t)=\frac{\partial Q(s,t)}{\partial s}$, we obtain the scaling relation
\[
g(s,t)
=
s^{2+\gamma-\frac{d}{\alpha}}\, g\!\left(1,\frac{t}{s}\right).
\]
We now prove that $f(x)>0$ for all $x\ge 1$. Integration by parts yields
\[
\int_0^1 w^{\gamma}(1-w)^{1-\frac{d}{\alpha}}\,dw
=
\frac{\gamma}{2-\frac{d}{\alpha}}
\int_0^1 w^{\gamma-1}(1-w)^{2-\frac{d}{\alpha}}\,dw,
\]
and therefore
\[
f(x)
=
\gamma^{-1}(x-1)^{2-\frac{d}{\alpha}}
+2\int_0^1 w^{\gamma-1}(1-w)^{2-\frac{d}{\alpha}}\,dw
-\int_0^1 w^{\gamma-1}(x+1-2w)^{2-\frac{d}{\alpha}}\,dw.
\]
In particular,
\[
f(1)
=
\big(2-2^{2-\frac{d}{\alpha}}\big)
\int_0^1 w^{\gamma-1}(1-w)^{2-\frac{d}{\alpha}}\,dw > 0,
\]
since $2-2^{2-\frac{d}{\alpha}}>0$.
To conclude, write
\[
f(x)
=
2\Bigg(\int_0^1 w^{\gamma-1}(1-w)^{2-\frac{d}{\alpha}}\,dw\Bigg)
\left[1-\frac{K(x)}{2}\right],
\]
where
\[
K(x)
=
\frac{
\int_0^1 w^{\gamma-1}(x+1-2w)^{2-\frac{d}{\alpha}}\,dw
-\gamma^{-1}(x-1)^{2-\frac{d}{\alpha}}
}{
\int_0^1 w^{\gamma-1}(1-w)^{2-\frac{d}{\alpha}}\,dw},
\quad x\ge1.
\]
A direct computation shows that $K'(x)<0$ for all $x\ge1$, since
$-1<1-\frac{d}{\alpha}<0$ and $(x+1-2w)>(x-1)$ for $0<w<1$.
Hence $K(x)\le K(1)<2$, which implies $f(x)>0$ for all $x\ge1$.
This completes the proof of Theorem \ref{Properties}.2.

\subsection{Proof of Theorem \ref{Properties}.3}
	
	Using integration by parts we obtain
	\begin{align*}
		\int_0^{s}w^{\gamma-1}(t-w)^{3-\da}\,dw &= 
		\frac{s^{\gamma}(t-s)^{3-\da}}{\gamma}
		+ \frac{\L(3-\da\R)}{\gamma}\int_0^sw^{\gamma}(t-w)^{2-\da}\,dw\\
		\int_0^{s}w^{\gamma-1}(s-w)^{3-\da}\,dw & = \frac{\L(3-\da\R)}{\gamma}\int_0^sw^{\gamma}(s-w)^{2-\da}\,dw,\\
		\int_0^s w^{\gamma-1}(t+s-2w)^{3-\da}\,dw &= \frac{s^{\gamma}(t-s)^{3-\da}}{\gamma}
		+ \frac{2\L(3-\da\R)}{\gamma}\int_0^s w^{\gamma}(t+s-2w)^{2-\da}\,dw.
	\end{align*}
	Hence, for $0\le s\le t$,
	\begin{align*}
		C(s,t)  = &\frac{3\L(3-\da\R)}{\gamma}\int_0^s w^{\gamma}(s-w)^{2-\da}\,dw
		+ 
		\frac{(t-s)^{3-\da}s^{\gamma}}{\gamma}
		+
		\frac{\L(3-\da\R)}{\gamma}\int_0^sw^{\gamma}(t-w)^{2-\da}\,dw \\
		& - \frac{\L(3-\da\R)}{\gamma}\int_0^sw^{\gamma}(s-w)^{2-\da}\,dw
		- \frac{s^{\gamma}(t-s)^{3-\da}}{\gamma}
		- \frac{2\L(3-\da\R)}{\gamma}\int_0^sw^{\gamma}(t+s-2w)^{2-\da}\,dw,
	\end{align*}
	and
	\begin{align*}
		C(s,s)=& \frac{4\L(3-\da\R)}{\gamma}\int_0^sw^{\gamma}(s-w)^{2-\da}\,dw
		- \frac{2^{3-\da}\L(3-\da\R)}{\gamma}\int_0^sw^{\gamma}(s-w)^{2-\da}\,dw\\
		=& \frac{\L(4-2^{3-\da}\R)\L(3-\da\R)}{\gamma}\int_0^sw^{\gamma}(s-w)^{2-\da}\,dw.
	\end{align*}
	Therefore,
	\begin{align} \nonumber
		\frac{\EE\L[\X(t)-\X(s)\R]^2}{\C} = &\ C(s,s)+C(t,t) -2C(s,t) \\ \nonumber
		= & \frac{\L(4-2^{3-\da}\R)\L(3-\da\R)}{\gamma}
		\L[
		\int_0^sw^{\gamma}(s-w)^{2-\da}\,dw + \int_0^tw^{\gamma}(t-w)^{2-\da}\,dw
		\R]
		\\ \nonumber
		& -\frac{6\L(3-\da\R)}{\gamma}\int_0^sw^{\gamma}(s-w)^{2-\da}\,dw - \frac{2s^\gamma\L(t-s\R)^{3-\da}}{\gamma}
		- \frac{2\L(3-\da\R)}{\gamma}\int_0^sw^{\gamma}(t-w)^{2-\da}\,dw \\ \nonumber
		&- \frac{2\L(3-\da\R)}{\gamma}
		\int_0^sw^{\gamma}\L[(t-w)^{2-\da}-(s-w)^{2-\da}\R]\,dw + \frac{2s^{\gamma}(t-s)^{3-\da}}{\gamma} \\ \nonumber
		& + \frac{4\L(3-\da\R)}{\gamma}
		\int_0^sw^{\gamma}(t+s-2w)^{2-\da}\,dw\\ \nonumber
		=& \L(-2-2^{3-\da}+2\R)\frac{\L(3-\da\R)}{\gamma}\int_0^sw^{\gamma}(s-w)^{2-\da}\,dw
\\ \nonumber &		+
		\L(4-2^{3-\da}-2-2\R)\frac{\L(3-\da\R)}{\gamma}\int_0^sw^{\gamma}(t-w)^{2-\da}\,dw \\ \nonumber
		& + \L(4-2^{3-\da}\R)\frac{\L(3-\da\R)}{\gamma}
		\int_s^tw^{\gamma}(t-w)^{2-\da}\,dw \\ \label{LabX} &
		+ 4\frac{\L(3-\da\R)}{\gamma} \int_0^sw^{\gamma}(t+s-2w)^{2-\da}\,dw.
	\end{align}
	Let $D_1=\L(4-2^{3-\da}\R)\frac{\L(3-\da\R)}{\gamma}$ and $D_2=\frac{\L(3-\da\R)}{\gamma}$. Notice that $D_1>0$ and $D_2>0$. By changing variables $w=(t-s)u+s$ we obtain
	\begin{align}\nonumber
		B_1 :=& \ D_1\int_s^t w^{\gamma}(t-w)^{2-\da}\,dw \\ \nonumber
		=&\ D_1\int_{0}^{1}\L[(t-s)u+s\R]^{\gamma}\L[t-s-(t-s)u\R]^{2-\da}(t-s)\,du \\ \label{X*}
		=&\ D_1(t-s)^{3-\da}\int_0^1\L[(t-s)u+s\R]^{\gamma}(1-u)^{2-\da}\,du \\ \nonumber
		>&\ D_1(t-s)^{3-\da}s^{\gamma}\int_{0}^{1}(1-u)^{2-\da}\,du =\frac{D_1}{3-\da}s^{\gamma}(t-s)^{3-\da}. 
	\end{align}
	For the remaining terms in \eqref{LabX} we have
	\begin{align*}\lefteqn{B_2:=
			2D_2\int_0^sw^{\gamma}\L[
			2(t+s-2w)^{2-\da}-2^{2-\da}(t-w)^{2-\da}-2^{2-\da}(s-w)^{2-\da}
			\R]\,dw} \\
		&\phantom{xx}=
		2D_2\int_0^s
		w^{\gamma}\L[
		2(t+s-2w)^{2-\da}-\L(2(t-w)\R)^{2-\da}-\L(2(s-w)\R)^{2-\da}
		\R]\,dw\\
		&\phantom{xxXx} \ge 0
	\end{align*}
	because the function $x\mapsto x^{2-\da}$ is concave. It follows that
	\begin{equation}\label{Estimacion1}
		\EE\L[\X(s)-\X(t)\R]^2\ \ge\ D_1\L(3-\da\R)^{-1}s^{\gamma}(t-s)^{3-\da}.
	\end{equation} 
	To obtain an upper bound for $\EE\L[\X(s)-\X(t)\R]^2$ we proceed as follows. From \eqref{X*},
	\begin{align}\label{Estimacion3}
		B_1 &\le D_1(t-s)^{3-\da}\int_0^1t^{\gamma}(1-u)^{2-\da}\,du = D_1\L(3-\da\R)^{-1}t^{\gamma}(t-s)^{3-\da}
		\le D_1\L(3-\da\R)^{-1}(t-s)^{3-\da}.
	\end{align}
	Since $w^{\gamma}\le s^{\gamma}$ for $0\le w\le s$,
	\begin{align}\nonumber
		B_2 &\le 2D_2 s^{\gamma}\int_0^s\L[
		2(t+s-2w)^{2-\da} - 2^{2-\da}(t-w)^{2-\da} - 2^{2-\da}(s-w)^{2-\da}\R]\,dw \\ \nonumber
		&= \frac{2D_2 s^{\gamma}}{3-\da}\L[
		(t+s)^{3-\da} -(t-s)^{3-\da} -2^{2-\da}t^{3-\da} +2^{2-\da}(t-s)^{3-\da}-2^{2-\da}s^{3-\da}
		\R]  \\ \nonumber
		& = \frac{2D_2 s^{\gamma}}{3-\da}\L[ 
		(t+s)^{3-\da} + \L(2^{2-\da}-1\R)(t-s)^{3-\da}-2^{2-\da}t^{3-\da} -2^{2-\da}s^{3-\da}
		\R] \\ \label{Estimacion2}
		& \le \frac{2D_2}{3-\da}\L(2^{2-\da}-1\R)(t-s)^{3-\da}
	\end{align}
	because the function $x\mapsto x^{3-\da}$ is convex due to   $1<3-\da<2$. Putting together \eqref{Estimacion3} and \eqref{Estimacion2} results in
	\begin{equation}\label{cotasuperior}
		\EE\L[\X(t)-\X(s)\R]^2 \le \L[
		D_1\L(3-\da\R)^{-1} + 2D_2 \L(3-\da\R)^{-1}\L(2^{2-\da}-1\R)
		\R](t-s)^{3-\da}.
	\end{equation}
	Finally, using \eqref{Estimacion1} and \eqref{cotasuperior} we obtain
	$$
	\frac{\L(4-2^{3-\da}\R)s^{\gamma}(t-s)^{3-\da}}{\gamma}\
	\le\ \EE\L[\X(t)-X(s)\R]^2\ \le\ \frac{2(t-s)^{3-\da}}{\gamma},
	$$
	and it follows from \cite[Corollary 2.1]{BGT} that the process $\{\X(t),\ t\ge0\}$ is not a semi-martingale on any interval $[\delta,1]$ for any $1>\delta>0$.
	
\subsection{Proof of Theorem \ref{Properties}.4}
We have proved above that 
$$
\EE\L[\X(t)-\X(s)\R]^2\ \le \ {\rm Const.}(t-s)^{3-\da},\quad 0\le s\le t.
$$
This shows that $\X$ has a.s. Hölder-continuous paths: for each $\varepsilon\in\L(0,\frac{3-\da}{2}\R)$ and each $L>0$ there exists a random variable  $K_{\varepsilon,L}\ge0$ such that
$$
|\X(t)-\X(s)|\ \le \ K_{\varepsilon,L}|t-s|^{\frac{3-\da}{2}-\varepsilon},\quad s,t,\in[0,L]\quad a.s.
$$

\subsection{Proof of Theorem \ref{Properties}.5}

Since $\X$ is Gaussian, its non-Markovianity follows from the fact that its covariance functional fails to satisfy the triangular identity
\[
C(s,u)\,C(t,t)=C(s,t)\,C(t,u),
\qquad 0\le s\le t\le u,
\]
which is a necessary and sufficient condition for Gaussian Markov processes.

\subsection{Proof of Theorem \ref{Properties}.6}

Since the covariance function \eqref{FormulaFinal} of $\X$ coincides with the covariance function $Q$ defined in (1.5) of \cite{ALEA}, Theorem 2.10(iv) in that paper applies with
\[
a=\gamma-1
\qquad\text{and}\qquad
b=2-\frac{d}{\alpha}.
\]
Observe that these parameters satisfy the required condition $b>0$.
Therefore, $\X$ exhibits long-range dependence. More precisely, for
$0\le r<v\le s<t$,
\[
\lim_{T\to\infty}
T^{d/\alpha}
\EE\big[(\X(t+T)-\X(s+T))(\X(v)-\X(r))\big]
=
\frac{2-d/\alpha}{\gamma(1+\gamma)}
(t-s)\big(v^{1+\gamma}-r^{1+\gamma}\big).
\]

\bigskip
{\bf \noindent Acknowledgements}\medskip

\noindent

The first and second authors were partially supported by grant CBF2023-2024-3562. They also acknowledge the hospitality of the Faculty of Mathematics and Informatics, Sofia University “St. Kliment Ohridski”, and the Institute of Mathematics and Informatics, Bulgarian Academy of Sciences, Sofia, Bulgaria, where this work was carried out.

The third author gratefully acknowledges the support provided by project UNITe BG16RFPR002-1.014-0004, funded by PRIDST, and partial support from funds allocated to Sofia University “St. Kliment Ohridski”, grant No. 80-10-41/2026.

\end{document}